\documentclass[12pt]{article}

\usepackage[english]{babel}
\usepackage[T1]{fontenc}
\usepackage[applemac]{inputenc}
\usepackage{lmodern}
\usepackage{microtype}
\usepackage{bbm}
\usepackage{tikz}
\usepackage{graphicx}
\usepackage{caption}
\usepackage{subcaption}
\usepackage{subfloat}
\usepackage{amsmath}
\usepackage{amssymb}
\usepackage{amsthm}
\usepackage{algorithm}
\usepackage[noend]{algpseudocode}
\usepackage{latexsym}
\usepackage{color}
\usepackage{dsfont}
\usepackage{appendix}
\newcommand{\argmin}{\operatornamewithlimits{argmin}}
\usepackage{hyperref}
\usepackage{float}
\usepackage{enumitem}
\usepackage{mathabx}
\usepackage{wrapfig}
\usepackage{graphicx} 
\graphicspath{{figures/}}

\DeclareSymbolFont{calletters}{OMS}{cmsy}{m}{n}
\DeclareSymbolFontAlphabet{\mathcal}{calletters}

\def\be{\begin{eqnarray}}
\def\ee{\end{eqnarray}}

\def\b*{\begin{eqnarray*}}
\def\e*{\end{eqnarray*}}

\def \L{\mathbb{L}}

\def \N{\mathbb{N}}
\def \P{{\mathbb P}}

\def \R{\mathbb{R}}

\def \a{\alpha}

\def \eps{\varepsilon}

\def\Ac{{\cal A}}
\def\Bc{{\cal B}}

\def\Dc{{\cal D}}

\def\Mc{{\cal M}}
\def\Nc{{\cal N}}

\def\Pc{{\cal P}}

\def\Wc{{\cal W}}
\def\Xc{{\cal X}}
\def\Yc{{\cal Y}}

\def\yb{{\bar y}}

\def \cl{{\rm cl\hspace{0.05cm}}}
\def \interior{{\rm int\hspace{0cm}}}

\def \conv{{\rm conv}}

\def \aff{{\rm aff}}

\def \argconc{{\rm argconc}}
\def \argmax{{\rm argmax}}

\def \Leb{{\L}}
\def \Ctn{{\rm{C}}}

\def \Ibf{{\mathbf I}}
\def \Sbf{{\mathbf S}}

\def \dist{{\rm dist}}

\def \zoom{{\rm zoom}}

\def\no{\noindent}
\def\x{\times}

\def\={\;=\;}
\def\.{\;.}

\def\eps{\varepsilon}

\def \1{{\bf 1}}
\def \ep{\hbox{ }\hfill{ ${\cal t}$~\hspace{-5.1mm}~${\cal u}$   } }

\def \proof{{\noindent \bf Proof. }}
\def \ep{\hbox{ }\hfill$\Box$}

 \def\normeL2#1{\left\|{#1}\right\|_{L^2}}

\title{Entropic approximation for multi-dimensional martingale optimal transport\thanks{The author gratefully acknowledges the financial support of the ERC 321111 Rofirm, and the Chairs Financial Risks (Risk Foundation, sponsored by Soci\'et\'e G\'en\'erale) and Finance and Sustainable Development (IEF sponsored by EDF and CA).}}

\author{Hadrien De March\thanks{CMAP, \'Ecole Polytechnique, hadrien.de-march@polytechnique.org.}}

\date{\today}

\begin{document}

\newtheorem{Theorem}{Theorem}[section]
\newtheorem{Lemma}[Theorem]{Lemma}
\newtheorem{Corollary}[Theorem]{Corollary}
\newtheorem{Proposition}[Theorem]{Proposition}
\newtheorem{Remark}[Theorem]{Remark}
\newtheorem{Example}[Theorem]{Example}
\newtheorem{Definition}[Theorem]{Definition}
\newtheorem{Assumption}[Theorem]{Assumption}

\maketitle


\abstract{We study the existing algorithms that solve the multidimensional martingale optimal transport. Then we provide a new algorithm based on entropic regularization and Newton's method. Then we provide theoretical convergence rate results and we check that this algorithm performs better through numerical experiments. We also give a simple way to deal with the absence of convex ordering among the marginals. Furthermore, we provide a new universal bound on the error linked to entropy.

\vspace{1cm}

}

\section{Introduction}

The problem of martingale optimal transport was introduced as the dual of the problem of robust (model-free) superhedging of exotic derivatives in financial mathematics, see Beiglb\"ock, Henry-Labord\`ere \& Penkner \cite{beiglbock2013model} in discrete time, and Galichon, Henry-Labord\`ere \& Touzi \cite{galichon2014stochastic} in continuous-time. This robust superhedging problem was introduced  by Hobson \cite{hobson1998robust}, and was addressing specific examples of exotic derivatives by means of corresponding solutions of the Skorokhod embedding problem, see \cite{cox2011robust,hobson2015robust,hobson2012robust}, and the survey \cite{hobson2011skorokhod}. 

Given two probability measures $\mu,\nu$ on $\R^d$, with finite first order moment, martingale optimal transport differs from standard optimal transport in that the set of all interpolating probability measures $\Pc(\mu,\nu)$ on the product space is reduced to the subset $\Mc(\mu,\nu)$ restricted by the martingale condition. We recall from Strassen \cite{strassen1965existence} that $\Mc(\mu,\nu)\neq\emptyset$ if and only if $\mu\preceq\nu$ in the convex order, i.e. $\mu(f)\le\nu(f)$ for all convex functions $f$. Notice that the inequality $\mu(f)\le\nu(f)$ is a direct consequence of the Jensen inequality, the reverse implication follows from the Hahn-Banach theorem.

This paper focuses on giving numerical aspects of martingale optimal transport for finite marginals. Henry-Labord\`ere \cite{henry2017model} used dual linear programming techniques to solve this problem, chosing well the cost functions so that the dual constraints were much easier to check. Alfonsi, Corbetta \& Jourdain noticed the difficulty, when going to higher dimension to get a discrete approximation of continuous marginals in convex order, that are still in convex order in higher dimension. So they mainly solve this problem, and then do several optimal transport resolutions with primal linear programming. Guo \& Ob{\l}{\'o}j \cite{guo2017computational} provide convergence results in the one dimensional setting of the discrete problem converges to the continuous problem, and they provide a Bregman projection scheme for solving the martingale optimal transport problem in the one dimensional setting. We also mention Tan \& Touzi \cite{tan2013optimal} who used a dynamic programming approach to solve a continuous-time version of martingale optimal transport.

The idea of using Bregman projection comes from classical optimal transport. Christian Leonard \cite{leonard2012schrodinger} was the first to have the idea of introducing an entropic penalization in an optimal transport problem. The entropic penalization makes this problem smooth and strictly convex and gives a Gibbs structure to the optimal probability, which has an explicit formula as a function of the dual optimizer. The unanimous adoption of entropic methods for solving optimal transport problems came from Marco Cuturi \cite{cuturi2013sinkhorn} who noticed that finding the dual solution of the entropic problem was equivalent to finding two diagonal matrices that made a full matrix bistochastic, therefore allowing to use the celebrated Sinkhorn algorithm.

Historically in classical optimal transport, the practitioners used linear programming algorithm to solve it, such as the Hungarian method \cite{kuhn1955hungarian}, the auction algorithm \cite{bertsekas1988auction}, the network simplex \cite{ahuja1988network}, we may also mention \cite{goldberg1989finding}. However, this method was so costly that only small problems could be treated because of the polynomial cost of linear programming algorithms. Later, Benhamou \& Brenier \cite{benamou2000computational} found another way of solving numerically the optimal transport problem by making it a dynamic programming problem with a final penalization on the mismatch of the final marginal of the dynamic process with the target marginal. For particular cases, it was also possible to use the Monge-Ampere equation. In the case of the square distance cost, Brenier \cite{brenier1991polar} proved that the optimal coupling is concentrated on a deterministic map, which was the gradient of a "potential" convex function $u$. When furthermore the marginals have densities with respect to the Lebesgue measure, we may prove that $u$ is a solution of the Monge-Ampere equation $\det D^2u = \frac{g\circ c_x(X,\cdot)^{-1}\circ\nabla u}{f}$, where $f$ is the density of $\mu$ and $g$ is the density of $\nu$. This equation satisfies a maximum principle, allowing to solve it in practice, see \cite{benamou2014numerical} and \cite{benamou2016monotone}. We also mention a smart strategy by Merigot \cite{merigot2011multiscale}, using semi-discrete transport. Levy \cite{levy2015numerical} introduced a Newton method to solve the semi-discrete problem very fast.

For the entropic resolution, Leonard \cite{leonard2012schrodinger} proved that the value of the entropic penalized optimal transport converged to the one of the unpenalized problem, while the optimal transports converged as well to a solution of the optimal transport. See \cite{carlier2017convergence} and \cite{cominetti1994asymptotic} for more precise studies of this convergence in particular cases. It have been observed by \cite{kosowsky1994invisible} that the entropic formulation was particularly useful for numerical resolution, as it allowed to use the celebrated Sinkhorn algorithm \cite{sinkhorn1967concerning}. The power of this technique has been rediscovered by \cite{cuturi2013sinkhorn}, and widely adopted by the community, see \cite{solomon2015convolutional}, \cite{rabin2015convex}, or \cite{thorpe2017transportation}. This method has already been adapted to different transport problems, such as Wasserstein barycenters \cite{agueh2011barycenters} and multi-marginal transport problems \cite{benamou2015iterative}, gradient flows problems \cite{peyre2015entropic}, unbalanced transport \cite{chizat2016scaling}, and one dimensional martingale optimal transport \cite{guo2017computational}.

The remarkable work by Schmitzer \cite{schmitzer2016stabilized} gives very practical considerations and tricks on how to actually make the Bregman projection algorithm converge fast and stay stable in practice. Cuturi \& Peyre \cite{cuturi2016smoothed} used a quasi-Newton method to solve the smooth entropic optimal transport. Their conclusion seems that the Sinkhorn algorithm is still more effective. However, \cite{brauer2017sinkhorn} use an inexact Newton method (i.e. including the use of the second derivative) and manage to beat the performance of the Sinkhorn algorithm. We also mention \cite{altschuler2017near} which introduces a "Greenkhorn algorithm" that outperforms the Sinkhorn algorithm according to their experiment, and similarly \cite{thibault2017overrelaxed} introduces an overrelaxed version of the Sinkhorn algorithm that squares the linear convergence coefficient, and accelerates the algorithm.

Our subsequent work differs from Guo \& Ob{\l}{\'o}j as we explain how to deal with higher dimension, give a more effective algorithm for martingale optimal transport by inexact Newton method. We also provide a speed of convergence for the Bregman projection algorithm, and explains how to deal with the lack of convex ordering of the marginals. Finally the universal bound that we give for the error linked to the entropy term is much sharper than the previous state-of-the-art. This bound may be extended to classical optimal transport for which it does not seem to be in the literature either.

In this paper we introduce several existing algorithms for solving martingale optimal transport such as linear programming, non-smooth semi-dual optimization, and Bregman projections. We introduce the smooth Newton algorithms, and the Newton semi-implied algorithm. Then we give some theoretical results on the speed of convergence of these algorithms, together with solutions to stabilize them and make them work in practice, like the preconditioning for the Newton method, or how to deal with marginals that are not in convex order. We provide new convergence rates for the entropic approximation of the martingale optimal transport, that are much better than the existing ones. The known result is an error of the order $\eps\big(\ln(N)-1\big)$, where $N$ is the size of the discretized grid, while we prove that we can get a result of order $\eps \frac{d}{2}$, where $d$ is the dimension of the space of the problem ($1$ or $2$ in this paper). These rates rely on very strong hypotheses that may be hard to check in practice. However we see on the numerical example that they are well verified in practice.

The paper is organized as follows. Section \ref{sect:preliminaries} gives the problem to solve, Section \ref{sect:algorithms} give the different algorithms that we will compare. In Section \ref{sect:practi}, we provide practical solutions to some usual problems, Section \ref{sect:conv_rate} provides theoretical convergence rates for the algorithms, Section \ref{sect:proofs} gathers the proofs of the theoretical results, and finally Section \ref{sect:numerics} contains numerical results.
\\
\newpage

\no {\bf Notation}\quad We fix an integer $d\ge 1$.

In all this paper, $\R^d$ is endowed with the Euclidean structure, the Euclidean norm of $x\in\R^d$ will be denoted $|x|$. Let $A\subset \R^d$ we denote $|A|$ the Lebesgue volume of $A$. The map $\iota_A$ is the map equal to $0$ on $A$, and $\infty$ otherwise. If $V$ is a topological affine space and $A\subset V$ is a subset of $V$, $\interior A$ is the interior of $A$, $\cl A$ is the closure of $A$, $\aff A$ is the smallest affine subspace of $V$ containing $A$, $\conv A$ is the convex hull of $A$, and $\dim(A):=\dim(\aff A)$. Let $(u_\eps)_{\eps>0},(v_\eps)_{\eps>0}\subset V$. We denote that $u_\eps = o(v_\eps)$ if $\lim_{\eps\to 0}\frac{|u_\eps|}{|v_\eps|}=0$. We further denote $u_\eps \ll v_\eps$. A classical property of $o(\cdot)$ is that
\be\label{eq:prop_o}
u_\eps = v_\eps+o(v_\eps)&\mbox{if and only if}&u_\eps = v_\eps+o(u_\eps).
\ee 

Let $x_0\in \R^d$, and $r>0$, we denote $\zoom_r^{x_0}:x\longmapsto x_0 + rx$, $B_r(x_0)$ is the closed ball centered in $x_0$ with radius $r$, and we only write $B_r$ when the center is $0$. Let $f:\R^d\longrightarrow \R$, we denote $\|f\|_\infty:=\sup_{x\in\R^d}f(x)$ its infinite norm, and $\|f\|^R_\infty:=\sup_{x\in B_R}f(x)$ its infinite norm when restricted to the ball $B_R$, for $R\ge 0$. Let $a,b\in\R^d$, we denote $a\otimes b:= ab^{t}=(a_ib_j)_{1\le i,j\le d}$, the only matrix in $\Mc_d(\R)$ such that for all $x\in\R^d$, we have $\left(a\otimes b\right) x = (b\cdot x)a$. Let $1\le k\le d+1$ and $u_1,...,u_k\in\R^d$, we denote $\det_\aff(u_1,...,u_k):=\left|\det\Big(\big(e_j\cdot(u_i-u_k)\big)_{1\le i,j\le k-1}\Big)\right|$, where $(e_j)_{1\le j\le k-1}$ is an orthonormal basis of $Vect\big(u_1-u_k,...,u_{k-1}-u_k\big)$.

Let $M\in\Mc_d(\R)$, a real matrix of size $d$, we denote $\det M$ the determinant of $M$. We also denote $Com(M)$ the comatrix of $M$: for $1\le i,j\le d$, $Com(M)_{i,j} = (-1)^{i+j}\det M^{i,j}$, where $M^{i,j}$ is the matrix of size $d-1$ obtained by removing the $i^{th}$ line and the $j^{th}$ row of $M$. Recall the useful comatrix formula:
\be\label{eq:comatrix}
Com(M)^t M = M Com(M)^t = (\det M) I_d.
\ee

As a consequence, whenever $M$ is invertible, $M^{-1} = \frac{1}{\det M}Com(M)^t$.

We denote $\Omega:=\R^d\times\R^d$ and define the two canonical maps
\b*
 X :(x,y)\in\Omega
 \longmapsto x\in\R^d
 &\mbox{and}&
 Y :(x,y)\in\Omega
 \longmapsto y\in\R^d.
 \e*
For $\varphi,\psi:\R^d\longrightarrow\bar\R$, and $h:\R^d\longrightarrow\R^d$, we denote 
\b*
\varphi\oplus\psi
:=
\varphi(X)+\psi(Y),
&\mbox{and}&
h^\otimes := h(X)\cdot(Y-X),
\e*
with the convention $\infty-\infty = \infty$.

For a Polish space $\Xc$, we denote by $\Pc(\Xc)$ the set of all probability measures on $\big(\Xc,\Bc(\Xc)\big)$. Let $\Yc$ be another Polish space, and $\P\in\Pc(\Xc\x\Yc)$. The corresponding conditional kernel $\P_x$ is defined by:
$$\P(dx,dy) = \P\circ X^{-1}(dx) \P_x(dy).$$
We also use this notation for finite measures. For a measure $m$ on $\Xc$, we denote $\Leb^1(\Xc,m):=\{f\in\Leb^0(\Xc):m[|f|]<\infty\}$. We also denote simply $\Leb^1(m):=\Leb^1(\bar\R,m)$.

\section{Preliminaries}
\label{sect:preliminaries}

Throughout this paper, we consider two probability measures $\mu$ and $\nu$ on $\R^d$ with finite first order moment, and $\mu \preceq \nu$ in the convex order, i.e. $\nu(f)\ge \mu(f)$ for all integrable convex $f$. We denote by $\Mc(\mu,\nu)$ the collection of all probability measures on $\R^d\times\R^d$ with marginals $\P\circ X^{-1}=\mu$ and $\P\circ Y^{-1}=\nu$. Notice that $\Mc(\mu,\nu)\neq\emptyset$ by Strassen \cite{strassen1965existence}.

For a derivative contract defined by a non-negative coupling function $c:\R^d\times\R^d\longrightarrow\R_+$, the martingale optimal transport problem is defined by:
 \b*
 \Sbf_{\mu,\nu}(c)
 &:=&
 \sup_{\P\in\Mc(\mu,\nu)}
 \P[c].
 \e*

The corresponding robust superhedging problem is
 \b*
 \Ibf_{\mu,\nu}(c)
 &:=&
 \inf_{(\varphi,\psi,h)\in\Dc_{\mu,\nu}(c)} \mu(\varphi)+\nu(\psi),
 \e*
where
 \b*
 \Dc_{\mu,\nu}(c)
 &:=&
 \big\{ (\varphi,\psi,h)\in\L^1(\mu)\x\L^1(\nu)\x\L^1(\mu,\R^d):
                                                           ~\varphi\oplus\psi+h^\otimes\ge c
 \big\}.~~~~~
 \e*
 The following inequality is immediate:
  \b*
\Sbf_{\mu,\nu}(c) \le \Ibf_{\mu,\nu}(c).
\e*
This inequality is the so-called weak duality. For upper semi-continuous coupling, we get from Beiglb\"ock, Henry-Labord\`ere, and Penckner \cite{beiglbock2013model}, and Zaev \cite{zaev2015monge} that there is strong duality, i.e. $\Sbf_{\mu,\nu}(c)= \Ibf_{\mu,\nu}(c)$. For any Borel coupling function bounded from below,  Beiglb\"ock, Nutz \& Touzi \cite{beiglbock2015complete} in dimension $1$, and De March \cite{de2018quasi} in higher dimension proved that duality holds for a quasi-sure formulation of dual problem and proved dual attainability thanks to the structure of martingale transports evidenced in \cite{de2017irreducible}.

Along all this paper, we assume that $\mu$ and $\nu$ are discrete, i.e. we may find finite $\Xc$ and $\Yc$ so that $\mu = \sum_{x\in\Xc}\mu_x\delta_x$, and $\nu = \sum_{y\in\Yc}\nu_y\delta_y$, so that all the coordinates of $\mu$ and $\nu$ are positive. Similarly, duality clearly holds thanks to the finiteness of the support, and the dual problem becomes discretized as well: for
$(\varphi,\psi,h)\in\Dc_{\mu,\nu}(c)$, we can denote $\varphi$, $\psi$, and $h$ as vectors $(\varphi(x))_{x\in\Xc}$, $(\varphi(y))_{y\in\Yc}$, and $(h_{i}(x))_{x\in\Xc,1\le i\le d}$.

To solve the martingale transport problem in practice, it seems necessary to discretize the problem. Guo \& Obl\'oj \cite{guo2017computational} prove that the martingale optimal transport problem with continuous $\mu$, $\nu$, and $c$ is a limit of this kind of discrete problem in dimension one under reasonable assumptions. This paper does not focus on proving the convergence of the discretized problem towards the continuous problem, we focus on how to solve the discretized problem.

\section{Algorithms}\label{sect:algorithms}

\subsection{Primal and dual simplex algorithm}

\subsubsection{Primal}

The natural strategy to solve this problem will be to use linear programming techniques such as simplex algorithm. One major problem with this approach is that the set $\Mc(\mu,\nu)$ may be empty, because in practice, the discretization of the marginals may break the convex ordering between then, thus making the set $\Mc(\mu,\nu)$ empty by Strassen theorem. This problem was relieved by Guo \& Ob{\l}{\'o}j \cite{guo2017computational}, and by Alfonsi, Corbetta \& Jourdain \cite{alfonsi2017sampling}. In \cite{guo2017computational}, they deal with the problem by replacing the convex ordering constraint by an approximate convex ordering constraint which is more resilient to perturbating the marginals. In \cite{alfonsi2017sampling}, they go beyond and gives several algorithms to find measures $\nu'$ (resp. $\mu'$) that are in convex order with $\mu$ (resp. with $\nu$) and satisfy some optimality criteria such as minimality of $\nu-\nu'$ (resp. $\mu-\mu'$) in terms of $p-$Wasserstein distance. We also give in Subsubsection \ref{subsubsect:convex_order} a technique to avoid this issue.

\subsubsection{Dual}

One huge weakness of the Primal algorithm is that the size of the problem is $|\Xc||\Yc|$, which is the size of $\Xc\x\Yc$, the support of the probabilities we consider. When $|\Xc|$ and $|\Yc|$ are big, it becomes a problem for memory storage. We notice that the number of constraints is $(d+1)|\Xc|+|\Yc|$, which is much smaller, because the dual functions $\varphi$, and $h$ are respectively in $\R^\Xc$ and in $(\R^d)^\Xc$, and the dual function $\psi$ lies in $\R^\Yc$. This is why in practice it makes sense to solve the Kuhn \& Tucker dual problem instead of the primal one. We will see considerations on the speed of convergence in Subsubsection \ref{subsubsect:converg_simplex}.

\subsection{Semi-dual non-smooth convex optimization approach}

It is well known from classical transport that solving directly the linear programming problem is too costly (see \cite{oberman2015efficient}) consequently, some alternative techniques have been developed like the Benamou-Brenier \cite{benamou2000computational} approach, which inspired Tan \& Touzi \cite{tan2013optimal} for the continuous time optimal transport problem. The idea consists in solving a Hamilton-Jacobi-Bellman problem with a penalization on the distance between the final marginal and $\nu$. Then an extension of this idea to our two-steps MOT problem gives the following resolution algorithm, suggested by Guo and Ob{\l}{\'o}j \cite{guo2017computational}. We denote $\Mc(\mu):=\{\P\in\Pc(\Omega):\P\circ X^{-1} = \mu,\mbox{ and }\P[Y|X] = X,\,\mu-\mbox{a.s.}\}$, and get
\b*
\Sbf_{\mu,\nu}(c)
 &:=&
 \sup_{\P\in\Mc(\mu,\nu)}
 \P[c]\\
 &=& \inf_{\psi\in\Leb^1(\nu) }\sup_{\P\in\Mc(\mu)}\P[c-\psi]+\nu[\psi]\\
 &=&\inf_{\psi\in\Leb^1(\nu) }\mu[(c(X,\cdot)-\psi)_{conc}(X)]+\nu[\psi]\\
 &=&\inf_{\psi\in\Leb^1(\nu) } V(\psi)
\e*

where $V(\psi) := \mu[(c(X,\cdot)-\psi)_{conc}(X)]+\nu[\psi]$ is a convex function in the variable $\psi$. Then the problem becomes a simple convex optimization problem. It seems appropritate in these conditions to solve the problem with using a classical gradient descent algorithm. It is proved in \cite{tan2013optimal} that $V$ has an explicit gradient. To give the explicit form of this gradient, we first need to introduce a notion of contact set. Let $f:\Yc\longmapsto \R$, as $\Yc$ is finite, $f_{conc}(x) = \inf_{f\le g\text{ affine}}g(x) = \sup\{\sum_i\lambda_i f(y_i):y_1,...,y_{d+1}\in\Yc,\lambda_1,...,\lambda_{d+1}\ge 0:\sum_i\lambda_iy_i = x\}$. By finiteness of $\Yc$, this supremum is a maximum. We denote $\argconc_f(x):=\argmax\{\sum_i\lambda_i f(y_i):y_1,...,y_{d+1}\in\Yc,\lambda_1,...,\lambda_{d+1}\ge 0:\sum_i\lambda_iy_i = x\}$. Then the subgradient of $V$ at $\psi$ is given by
$$\partial V(\psi) = \left\{\sum_{x\in\Xc}\mu_x\sum_i\lambda_i(x)\delta_{y_i(x)}-\nu:\big(y(x),\lambda(x)\big)\in \argconc_{c(x,\cdot)-\psi}(x),\mbox{ for all }x\in\Xc\right\}$$
Notice that this set is a singleton for a.e. $\psi\in\Leb^1(\Yc)$, as $V$ is a convex function in finite dimensions. Then with high probability, on each gradient step, the function $V$ will be differentiable on this point. In practice there is always uniqueness after the first step.

\subsection{Entropic algorithms}

\subsubsection{The entropic problem in optimal transport}

In practice, this problem is added some regularity by the addition of an entropic penalization (see Leonard \cite{leonard2012schrodinger}, Cuturi \cite{cuturi2013sinkhorn}). Let $\eps>0$,
\b*
\Sbf_{\mu,\nu}^\eps(c):= \sup_{\P\in\Pc(\mu,\nu)}\P[c]-\eps H(\P|m_0),
\e*
where $H(\P|m_0):= \int_\Omega \left(\ln\left(\frac{d\P}{dm_0}\right)-1\right)\frac{d\P}{dm_0}m_0(d\omega)$. The measure $m_0$ is the "reference measure", we assume that it may be decomposed as $m_0 := m^{\Xc}\otimes m^{\Yc}$ such that $\mu$ is dominated by $m^{\Xc}\in\Mc(\R^d)$, and $\nu$ is dominated by $m^{\Yc}\in\Mc(\R^d)$.
For this text we chose $m_0:=\sum_{(x,y)\in\Xc\x\Yc}\delta_{(x,y)}$. By the finiteness of the supports of $\mu$ and $\nu$, we know that $\P$ is absolutely continuous with respect to $m_0$. Denote $p:=\frac{d\P}{d m_0}$, and abuse notation writting $p\in\Pc(\mu,\nu)$. Then $H(\P| m_0):= \sum_{(x,y)\in\Xc\x\Yc}\left(\ln p(x,y)-1\right)p(x,y)$. This problem can be written with Lagrange multipliers,
\b*
\sup_{\P\in\Pc(\mu,\nu)}\P[c]-\eps H(\P| m_0) = \inf_{(\varphi,\psi)\in\Leb^1(\mu)\x\Leb^1(\nu)}\sup_{\P\in\Pc(\mu,\nu)}\P[c-\varphi\oplus\psi]-\eps H(\P| m_0)+\mu[\varphi]+\nu[\psi].
\e*
which leads to an explicit Gibbs form for the optimal kernel $p$. Then as the supports are finite we easily get the shape of the optimizer $p(x,y)=\exp\left(-\frac{\varphi(x)+\psi(y)-c(x,y)}{\eps}\right)$, and the associated dual problem becomes
\b*
\Ibf_{\mu,\nu}^\eps
(c):= \inf_{(\varphi,\psi)\in\Leb^1(\mu)\x\Leb^1(\nu)}\mu[\varphi]+\nu[\psi]+\eps\sum_{x,y}\exp\left(-\frac{\varphi(x)+\psi(y)-c(x,y)}{\eps}\right).
\e*

One important property that we need is the $\Gamma-$convergence. We say that $F_\eps$ $\Gamma-$converges to $F$ when $\eps \longrightarrow 0$ if for all sequence $\eps_n\longrightarrow 0$, we have

\no{\rm (i)} For all sequences $x_n\longrightarrow x$, we have $F(x)\ge\limsup_n F_{\eps_n}(x_n)$.

\no{\rm (ii)} There exists a sequence $x_n\longrightarrow x$ such that $F(x)\le\liminf_n F_{\eps_n}(x_n)$.

The $\Gamma-$convergence implies that $\min F_n\longrightarrow F$, when $n\longrightarrow\infty$, and that if $x_n$ is a minimizer of $F_n$ for all $n\ge 1$, and if $x_n\longrightarrow x$, then $x$ is a minimizer of $F$. Leonard \cite{leonard2012schrodinger} proved this $\Gamma-$convergence of the penalized problem to the optimal transport problem.

\subsubsection{The Bregman iterations algorithm}

Coupled with the Sinkhorn algorithm \cite{sinkhorn1967concerning} introduced by Marco Cuturi for optimal transport \cite{cuturi2013sinkhorn}, this method allows an exponentially fast approximated resolution. Notice that the operator $V_\eps(\varphi,\psi):=\mu[\varphi]+\nu[\psi]+\eps\sum_{x,y}\exp\left(-\frac{\varphi(x)+\psi(y)-c(x,y)}{\eps}\right)$ is smooth convex. The Euler-Lagrange equations $\partial_\varphi V_\eps = 0$ (resp. $\partial_\psi V_\eps = 0$) are exactly equivalent to the marginal relations $\P\circ X^{-1} = \mu$ (resp. $\P\circ Y^{-1} = \nu$). It was noticed in \cite{cuturi2013sinkhorn} that these partial optimizations can be obtained in closed form:
\b*
\varphi(x) = \eps\ln\left(\frac{1}{\mu_x}\sum_{y}\exp\left(-\frac{\psi(y)-c(x,y)}{\eps}\right)\right),
\e* 
and
\b* 
\psi(y) = \eps\ln\left(\frac{1}{\nu_y}\sum_{x}\exp\left(-\frac{\varphi(x)-c(x,y)}{\eps}\right)\right).
\e*
By iterating these partial optimization, we obtain the so-called Sinkhorn algorithm (see \cite{sinkhorn1967concerning}) that is equivalent to a block optimization of the smooth function $V_\eps$ which dual is called Bregman projection \cite{bregman1967relaxation}, and converges exponentially fast, see Knight \cite{knight2008sinkhorn}.

\subsubsection{The entropic approach for the one period martingale optimal transport problem}

As observed by Guo \& Ob{\l}{\'o}j \cite{guo2017computational} in dimension 1, the Sinkhorn algorithm can be extended to the martingale optimal transport problem. With exactly the same computations, we get
\b*
\Sbf_{\mu,\nu}^\eps(c)&:=& \sup_{\P\in\Mc(\mu,\nu)}\P[c]-\eps H(\P| m_0)\\
=\Ibf_{\mu,\nu}^\eps(c)&:=& \inf_{(\varphi,\psi,h)\in\Leb^1(\mu)\x\Leb^1(\nu)\x\Leb^0(\R^d)}\mu[\varphi]+\nu[\psi]\\
&&+\eps\sum_{x,y}\exp\left(-\frac{\varphi(x)+\psi(y)+h^\otimes(x,y)-c(x,y)}{\eps}\right).
\e*
First notice that the $\Gamma-$convergence still holds in this easy finite case.
\begin{Proposition}
Let $F_\eps : \P\in\Mc(\mu,\nu)\longmapsto \P[c] - \eps H(\P|m_0)$. For $\eps>0$, $F_\eps$ is strictly concave upper semi-continuous. Furthermore, $F_\eps$ $\Gamma-$converges to $F_0$ when $\eps\longmapsto 0$.
\end{Proposition}
\proof
This $\Gamma-$convergence is easy by finiteness as the entropy is bounded by $\ln(|\Xc||\Yc|)-1$ when $\P$ is a probability measure.
\ep\\

We denote $\Delta:=\varphi\oplus\psi+h^\otimes-c$, the convex function to minimize becomes
$$V_\eps(\varphi,\psi,h):=\mu[\varphi]+\nu[\psi]+\eps\sum_{x,y}\exp\left(-\frac{\Delta(x,y)}{\eps}\right).$$
Then the Sinkhorn algorithm is complemented by another step so as to account for the martingale relation:
\be\label{eq:partial_opt}
\varphi(x) &=& \eps\ln\left(\frac{1}{\mu_x}\sum_{y}\exp\left(-\frac{\psi(y)+h(x)\cdot(y-x)-c(x,y)}{\eps}\right)\right),\\
\psi(y) &=& \eps\ln\left(\frac{1}{\nu_y}\sum_{x}\exp\left(-\frac{\varphi(x)+h(x)\cdot(y-x)-c(x,y)}{\eps}\right)\right),\nonumber\\
0 &=& \frac{1}{\mu_x}\sum_{y}\exp\left(-\frac{\Delta(x,y)}{\eps}\right)(y-x)\nonumber\\
&=&\frac{-1}{\mu_x}\frac{\partial}{\partial h(x)}\eps\sum_{y}\exp\left(-\frac{\Delta(x,y)}{\eps}\right).\nonumber
\ee
Notice that the martingale step is not closed form and is only implied. However, it may be computed almost as fast as $\varphi$, and $\psi$, thanks to the Newton algorithm applied to each smooth strongly convex function $F_x$ of $d$ variables given, for each $x\in\Xc$ with its derivatives by
\be\label{eq:func_martingale}
F_x(h)&:=&\eps\sum_{y}\exp\left(-\frac{\Delta(x,y)}{\eps}\right),\\
\nabla F_x(h)&=&-\sum_{y}\exp\left(-\frac{\Delta(x,y)}{\eps}\right)(y-x),\nonumber\\
D^2 F_x(h)&=&\frac{1}{\eps}\sum_{y}\exp\left(-\frac{\Delta(x,y)}{\eps}\right)(y-x)\otimes(y-x).\nonumber
\ee
Notice that the optimization of $F_{x_1}$ and $F_{x_2}$ are independent for $x_1\neq x_2$.

\subsubsection{Truncated Newton method}

For these problems, it may make sense to use a Newton method, as the problems are smooth, and the Newton method converges very fast. For very highly dimensional problems (here $(d+1)|\Xc|+|\Yc|$), the inversion of the hessian is too much costly. Then it is in general preferred to use quasi-Newton. Instead of computing the Newton step $D^2V_\eps^{-1}\nabla V_\eps$, we use a conjugate gradient algorithm to find by iterations a vector $p\in\Dc_{\Xc,\Yc}$ such that $|D^2V_\eps p - \nabla V_\eps|$ is small enough, generally in practice this quantity is chosen to be smaller than $\min\left(\frac12,\sqrt{|\nabla V_\eps|}\right)$.

The conjugate gradient algorithm approximates the solution of the equation $Ax = b$ by solving it "direction by direction" along the most important direction, until a stopping criterion is reached. The exact algorithm may be found in \cite{wright1999numerical}.

\subsubsection{Implied truncated Newton method}

Some instabilities may appear from Newton steps as any term of the form $\exp(X/\eps)$ can easily explode when $\eps$ is very small and $X>0$. The dimension may also make the conjugate gradient from the quasi-Newton algorithm slow. A good way to avoid this problem and exploit the near-closed formula for the optimal $\varphi$ and $h$ when $\psi$ is fixed, or optimal $\psi$ when $\varphi$ and $h$ are fixed.

Instead of applying the truncated method to $V_\eps(\varphi,\psi,h)$, we apply the truncated Newton method to $\widetilde{V}_\eps(\psi) := \min_{\varphi,h}V_\eps(\varphi,\psi,h)$. It is elementary that with these definitions we have
$$\inf_{\varphi,\psi,h}V_\eps(\varphi,\psi,h) = \inf_{\psi}\widetilde{V}_\eps(\psi).$$
Doing this variable implicitation is easy by the fact that we have a closed formula for $\varphi$ and a quasi-closed formula for $h$. It brings the great advantage of having the first marginal and the martingale relationship verified, this fact will be exploited in Subsubsection \ref{subsubsect:convex_order}.

Now we give a general framework that allows to use variables implicitation. The following framework should be used with $F = V_\eps$, $x = \psi$, and $y = (\varphi,h)$. Proposition \ref{prop:varimpl} below provides the appropriate convexity result together with the closed formulas for the two first derivatives of $\widetilde{V}_\eps$ that are necessary to apply the truncated Newton algorithm. Let $\Ac$ and $\Bc$ finite dimensional spaces and $F:\Ac\x\Bc\longrightarrow \R$, we say that $F$ is $\alpha-$convex if
\b*
\lambda F(\omega_1)+(1-\lambda)F(\omega_2)-F\big(\lambda \omega_1+(1-\lambda)\omega_2\big)\ge \alpha\frac{\lambda(1-\lambda)}{2}|\omega_1-\omega_2|^2,
\e*
for all $\omega_1,\omega_2\in\Ac\x\Bc$, and $0\le \lambda\le 1$. The case $\alpha = 0$ corresponds to the standard notion of convexity.
\begin{Proposition}\label{prop:varimpl}
Let $F:\Ac\x\Bc\longrightarrow \R$ be a $\alpha-$convex function. Then the map $\tilde{F}:x\longmapsto\inf_{y\in\Bc}F(x,y)$ is $\alpha-$convex. Furthermore, if $\alpha>0$ and $F$ is $\Ctn^2$, then $y(x):=\argmin_y F(x,y)$ is unique and we have
\b*
\nabla y(x) &=& -\partial^2_{y}F^{-1}\partial^2_{yx}F\big(x,y(x)\big),\\
\nabla \tilde{F}(x) &=& \partial_x F\big(x,y(x)\big),\\
D^2 \tilde{F}(x) &=& \left(\partial^2_xF-\partial^2_{xy}F\partial^2_{y}F^{-1}\partial^2_{yx}F\right)\big(x,y(x)\big).
\e*
\end{Proposition}

The proof of Proposition \ref{prop:varimpl} is reported in Subsection \ref{subsect:conv_min}.

\begin{Remark}
The matrix $\partial^2_{xy}F\partial^2_{y}F^{-1}\partial^2_{yx}F$ is symmetric positive definite, therefore the curvature of the function $\tilde{F}$ is reduced by the implicitation process, making heuristically the minimization easier. This fact is also observed in practice.
\end{Remark}

This method shall be used for the optimization of $V_\eps$, but also for the optimization of $F_x$ that gives the martingale step, see \eqref{eq:func_martingale}. Indeed the value of $\varphi(x)$ does not change the martingale optimality of $F_x$. We provide these important formulas.

\paragraph{The map $\widetilde V_\eps$ and its derivatives:}

Let $\psi\in \R^\Yc$, we denote $(\widehat\varphi^\eps_\psi,\widehat{h}^\eps_\psi):= \argmin_{\varphi,h}V_\eps(\varphi,\psi,h)$, that are unique and may be found in quasi-closed form from \eqref{eq:partial_opt}. Now we give the formula for $\widetilde{V}_\eps$ and its derivatives. We directly get from Proposition \ref{prop:varimpl} that
\b*
\widetilde{V}_\eps(\psi) &= &V_\eps\left(\widehat{\varphi}^\eps_\psi,\psi,\widehat{h}^\eps_\psi\right),\\
\nabla\widetilde{V}_\eps(\psi) &=& \partial_\psi V_\eps\left(\widehat{\varphi}^\eps_\psi,\psi,\widehat{h}^\eps_\psi\right),\\
D^2\widetilde V_\eps &=&\left(\partial^2_\psi V_\eps-\partial^2_{\psi,\varphi}V_\eps(\partial^2_{\varphi}V_\eps)^{-1}\partial^2_{\varphi,\psi}V_\eps-\sum_{i=1}^d\partial^2_{\psi,h_i}V_\eps(\partial^2_{h_i}V_\eps)^{-1}\partial^2_{h_i,\psi}V_\eps\right)\left(\widehat{\varphi}^\eps_\psi,\psi,\widehat{h}^\eps_\psi\right).
\e*
The last additive decomposition of $\partial_{(\varphi,h)^2} V_\eps^{-1}$ stems from the fact that $\partial_{(\varphi,h)^2}V_\eps\left(\widehat{\varphi}^\eps_\psi,\psi,\widehat{h}^\eps_\psi\right)$ is diagonal. Indeed, $V_\eps$ is a sum of functions of $\left(\varphi(x),h(x)\right)$ for $x\in\Xc$, and the crossed derivative $\partial_{\varphi(x),h(x)}V_\eps = \sum_y (y-x)\exp\left(-\frac{\Delta(x,y)}{\eps}\right)$ cancels at $\left(\widehat{\varphi}^\eps_\psi(x),\widehat{h}^\eps_\psi(x)\right)$ by the martingale property induced by the optimality in $h(x)$. The same holds for $\partial_{h_i(x),h_j(x)}V_\eps$ for $i\neq j$. We denote $\Delta_\psi^\eps:=\widehat\varphi^\eps_\psi\oplus\psi+\left(\widehat h_\psi^\eps\right)^\otimes-c$, and we have
\b*
\widetilde{V}_\eps(\psi) &=&\mu\left[\widehat\varphi^\eps_\psi\right]+\nu[\psi]+\eps\sum_{x,y}\exp\left(-\frac{\Delta^\eps_\psi(x,y)}{\eps}\right),\\
\nabla\widetilde{V}_\eps(\psi) &=&\left(\nu_y -\sum_{x}\exp\left(-\frac{\Delta^\eps_\psi(x,y)}{\eps}\right)\right)_{y\in\Yc},\\
D^2\widetilde V_\eps &=&\eps^{-1}diag\left(\sum_x\exp\left(-\frac{\Delta^\eps_\psi(x,y)}{\eps}\right)\right)\\
&&-\eps^{-1}\sum_x
\left(\sum_y\exp\left(-\frac{\Delta^\eps_\psi(x,y)}{\eps}\right)\right)^{-1}\\
&&\quad\quad\x\Bigg(\exp\left(-\frac{\Delta^\eps_\psi(x,y_1)}{\eps}\right)\exp\left(-\frac{\Delta^\eps_\psi(x,y_2)}{\eps}\right)\Bigg)_{y_1,y_2\in\Yc}\\
&&-\eps^{-1}\sum_x\sum_{i=1}^{d}\left(\sum_y(y_i-x_i)^2\exp\left(-\frac{\Delta^\eps_\psi(x,y)}{\eps}\right)\right)^{-1}\\
&&\quad\quad\x\Bigg((y_1-x)_i\exp\left(-\frac{\Delta^\eps_\psi(x,y_1)}{\eps}\right)(y_2-x)_i\exp\left(-\frac{\Delta^\eps_\psi(x,y_2)}{\eps}\right)\Bigg)_{y_1,y_2\in\Yc}.
\e*
Notice that for the conjugate gradient algorithm, we only need to be able to compute the product $\left(D^2\widetilde V_\eps\right) p$ for $p\in\R^\Yc$. Then if the RAM is not sufficient to store the whole matrix $D^2\widetilde V_\eps$, it may be convenient to only store $D_\psi := diag\left(\sum_x\exp\left(-\frac{\Delta^\eps_\psi(x,y)}{\eps}\right)\right)$, $D_\varphi^{-1} := diag\left(\sum_y\exp\left(-\frac{\Delta^\eps_\psi(x,y)}{\eps}\right)\right)^{-1}$, and $D_{h_i}^{-1} := diag\left(\sum_y(y_i-x_i)^2\exp\left(-\frac{\Delta^\eps_\psi(x,y)}{\eps}\right)\right)^{-1}$ for all $1\le i \le d$. Then we compute $\left(D^2\widetilde V_\eps\right) p$ in the following way:
\b*
p_\psi &:=&\left(D_\psi\right) p\\
p_\varphi &:=& \left(\sum_y \exp\left(-\frac{\Delta^\eps_\psi(x,y)}{\eps}\right)p_y\right)_{x\in\Xc},\\
p_{h_i}&:=&\left(\sum_y (y_i-x_i)\exp\left(-\frac{\Delta^\eps_\psi(x,y)}{\eps}\right)p_y\right)_{x\in\Xc},\\
p_\varphi' &:=& \left(D_\phi^{-1}\right)p_\varphi,\\
p_{h_i}'&:=&\left(D_{h_i}^{-1}\right)p_{h_i},\\
p_\varphi'' &:=& \left(\sum_x \exp\left(-\frac{\Delta^\eps_\psi(x,y)}{\eps}\right)(p_\varphi'')_x\right)_{y\in\Yc},\\
p_{h_i}''&:=&\left(\sum_x (y_i-x_i)\exp\left(-\frac{\Delta^\eps_\psi(x,y)}{\eps}\right)(p_{h_i}'')_x\right)_{y\in\Yc},\\
\left(D^2\widetilde V_\eps\right) p&=&\eps^{-1}\left(p_\psi - p_\varphi''-\sum_{i=1}^dp_{h_i}''\right).
\e*

\paragraph{The map $\widetilde F_x$ and its derivatives:}

In this paragraph we fix $\psi\in\R^\Yc$ and $\eps>0$. Recall the map $F_x$ from \eqref{eq:func_martingale}. This map may be seen as a function of $\big(\varphi(x),h(x)\big)$. Then we set
\b*
\widetilde F_x(h):=\min_{\varphi(x)\in\R}\mu_x\varphi(x)+\eps\sum_y\exp\left(-\frac{\varphi(x)+\psi+h\cdot(y-x)-c(x,y)}{\eps}\right).
\e*
The optimizer is given by \eqref{eq:partial_opt}, hence by the closed formula
\b*
\widehat\varphi_h(x):=\argmin_{\varphi(x)\in\R} \mu_x\varphi(x)+F_x\big(\varphi(x),h\big) = \eps\ln\left(\frac{1}{\mu_x}\sum_{y}\exp\left(-\frac{\psi(y)+h\cdot(y-x)-c(x,y)}{\eps}\right)\right).
\e*
A direct application of Proposition \ref{prop:varimpl} gives
\b*
\widetilde F_x(h)&=&\min_{\varphi(x)\in\R}\mu_x\widehat\varphi_h(x)+\eps\sum_y\exp\left(-\frac{\widehat\Delta_h(x,y)}{\eps}\right),\\
\nabla\widetilde F_x(h)&=&-\sum_y(y-x)\exp\left(-\frac{\widehat\Delta_h(x,y)}{\eps}\right),\\
D^2\widetilde F_x(h)&=&\eps^{-1}\left(\sum_y(y-x)^{2\otimes}\exp\left(-\frac{\widehat\Delta_h(x,y)}{\eps}\right)
-\mu_x^{-1}\left(\sum_y(y-x)\exp\left(-\frac{\widehat\Delta_h(x,y)}{\eps}\right)\right)^{2\otimes}\right),
\e*
where we denote $\widehat\Delta_h(x,y):=\widehat\varphi_h(x)+\psi(y)+h\cdot(y-x)-c(x,y)$, and $u^{2\otimes}:=u\otimes u$ for $u\in\R^d$.

\section{Solutions to practical problems}\label{sect:practi}

\subsection{Preventing numerical explosion of the exponentials}

As we want to make $\eps$ go to $0$, all the terms like $\exp\left(\frac{\cdot}{\eps}\right)$ tend to explode numerically. Here are the different risks that we have to deal with, and how we deal with them.

First the Newton algorithm is very local, and nothing guarantees that after one iteration, the value function will not explode. From our practical experience, the algorithm tends to explode for $\eps < 10^{-3}$. Notice that the numerical experiment given by \cite{brauer2017sinkhorn} does not go beyond $10^{-3}$, we may imagine that this is because they do not use the variable implicitation technique. Furthermore, we notice from our numerical experimenting that this variable implicitation, additionaly to the stabilizing the numerical scheme, makes the convergence of the Newton algorithm much faster. Moreover, impliciting in $\varphi$ and $h$ is much more effective than impliciting in $\psi$, even though we have to do the implicitation in $h$ which is much more costly than the implicitation in $\psi$.

For the computation of the implicitations \eqref{eq:partial_opt}, the computation of the formula $\varphi(x) = \eps\ln\left(\frac{1}{\mu_x}\sum_{y}\exp\left(-\frac{\psi(y)+h(x)\cdot(y-x)-c(x,y)}{\eps}\right)\right)$ should be done as follows to prevent numerical explosion. First we compute $M_x := \max_{y\in\Yc}\left\{-\frac{\psi(y)+h(x)\cdot(y-x)-c(x,y)}{\eps}\right\}$, and then the computation that we do effectively is
\be\label{eq:exp_max}
\varphi(x) &=& M_x+\eps\ln\left(\sum_{y}\exp\left(-\frac{\psi(y)+h(x)\cdot(y-x)-c(x,y)}{\eps}-M_x\right)\right)-\eps \ln \mu_x
\ee
In \eqref{eq:exp_max}, the exponential arguments are always smaller than $1$, and one of them is equal to $1$, then any explosion makes the exponential be totally negligible when compared to $\exp(0)=1$, this computation rule makes it very stable. Notice also the separation of $\ln \mu_x$ that allows to treat the case when the value of $\mu_x$ is extremely low (like for exemple when you discretise a Gaussian measure on a grid) even if in this case, it may be smarter to just remove the value from the grid.

Notice that the variable implicitation should also be used during each partial optimisation in $h(x)$ for $x\in\Xc$, as this Newton algorithm is highly susceptible to explode as well. The implicitation simply consists in minimizing in $\varphi(x)$ the maps $F_x$ from \eqref{eq:func_martingale}, and has a closed form.

Another thing to take care of about $h$ is the initial value taken for the next partial optimization of $V_\eps$ in $h$. On a first hand, chosing the last value for $h$ helps diminishing the number of steps for the optimization. Also, when $\eps$ is very small, even with the implicitation, the Newton optimization may get hard if the initial value is too far from the optimum. 

\subsection{Customization of the Newton method}

\subsubsection{Preconditioning}

The conjugate gradient algorithm used to compute the search direction for the Newton algorithm has a convergence rate given by $|x_k-x^*|_A\le 2\left(\frac{\sqrt{\kappa(A)}-1}{\sqrt{\kappa(A)}+1}\right)^k|x_0-x^*|_A$, where $x_k$ is the $k-$th iterate, $x^*$ is the solution of the problem, $|x|_A := x^tAx$ is the Euclidean norm associated to the scalar product $A$, and $\kappa(A):=\Vert A\Vert\Vert A^{-1}\Vert$ is the conditioning of $A$. This conditioning is the fundamental parameter for this convergence speed. When $\eps$ is getting small, the conditioning raises. We also observe on the numerics that is happens when the marginals have a thin tail (e.g. Gaussian distributions).
The simplest way of dealing with this conditioning problem consists in applying a "preconditioning" algorithm. We find a matrix $P$ that is easy to invert (for example a diagonal matrix) and we use the fact that solving $Ax = b$ is equivalent with solving $P^tAPx' = P^tb$, where $x' := P^{-1}x$. We use the most classical and simple preconditioning which consists in taking $P:= \sqrt{diag(A)}^{-1}$. See \cite{wright1999numerical} for the precise algorithm.


\subsubsection{Line search}

An important advantage of the Bregman projection algorithm over the primitive Newton algorithm is that $V_\eps$ is a Lyapunov function as the steps only consist of block minimizations of this function, whereas the Newton step may get very wrong and lose the optimal region if we are not close enough to the minimum. However in practice, some ingredients need to be added to the Newton step. Indeed, once the direction of search is decided by the conjugate gradient algorithm, in practice it is necessary to make a line search algorithm, i.e. to find a point on the line on which the value function $V_\eps$ is strictly smaller, and so does the directional gradient absolute value $|\nabla V_\eps\cdot p|$, where $p$ is the descent direction. This "descent" condition is called the Wolfe condition. A very good line search algorithm that is commonly used in practice is detailed in \cite{wright1999numerical}.

\begin{Remark}
Notice that if a value is rejected by the line search, it is important to throw away the value of $h$ given by this wrong point, and to come back to the last value of $h$ corresponding to a point that was not rejected by the line search.
\end{Remark}

\subsection{Penalization}

\subsubsection{The penalized problem}

The dual solution may not be unique, which may lead to numerical unstabilities. As an example we may add any constant to $\varphi$ while subtracting it to $\psi$ without changing the value of $V$. A straightforward solution is to add a penalization to the minimization problem. I.e. we have the new problem
\be\label{pb:penoptimization}
\min_{\psi\in\R^\Yc} \widetilde{V}_\eps(\psi)+\a f(\psi)
\ee
where $f$ is a strictly convex superlinear function, so that there is a unique minimum by the fact that the gradient of $\widetilde{V}_\eps$ is a difference of probabilities, which proves that this convex function is Lipschitz, whence the strict convexity and super-linearity of $\widetilde{V}_\eps(\psi)+\a f(\psi)$. In practice we take $f(\psi):=\frac12\sum_{y\in\Yc}a_y \psi(y)^2$, for some $a\in\R^\Yc$, so that $\nabla f(\psi) = \sum_{y\in\Yc}a_y\psi(y)e_y$, where $(e_y)_{y\in\Yc}$ is the canonical basis, and $D f(\psi) = {\rm diag}(a)$ have these easy closed expressions. In practice we take $a = (1)$, $a = \nu$, $a = \nu^2$, or $a = \nu/\psi_0$, where $\varphi_0$ is a fixed estimate of $\psi$ from the last step of $\eps-$scaling (see Subsection \ref{subsect:eps_scaling}).

\subsubsection{Marginals not in convex order}\label{subsubsect:convex_order}

Problem \ref{pb:penoptimization} allows to solve the problem of mismatch in the convex ordering thanks to the following theorem that allows for probability measures $\mu,\nu$ not in convex order to find another probability measure $\widetilde{\nu}$ in convex order with $\mu$ that satisfies some optimality criterion, for example in terms of distance from $\nu$.

\begin{Theorem}\label{thm:marginal_limit}
Let $(\mu,\nu)\in \Pc(\Xc)\x\Pc(\Yc)$ not in convex order. Let $\nu_\alpha:=\P_\alpha\circ Y^{-1}$, where $\P_\alpha$ is the optimal probability for Problem \eqref{pb:penoptimization}, where $f$ is a super-linear, differentiable, strictly convex, and $p-$homogeneous function $\R^\Yc\longrightarrow\R$ for some $p>1$. Then $\nu_\alpha \longrightarrow \nu_l$ when $\alpha\longrightarrow 0$, for some $\nu_l\succeq_c\mu$ satisfying
\b*
f^*(\nu_l-\nu) = \min_{\tilde{\nu}\succeq_c\mu}f^*(\tilde{\nu}-\nu).
\e*
\end{Theorem}

The proof of Theorem \ref{thm:marginal_limit} is reported in Subsection \ref{subsect:marginal_limit}. Notice that for $f(\psi):=\frac12\sum_{y\in\Yc}a_y \psi(y)^2$, we have $f^*(\gamma)=\frac12\sum_{y\in\Yc}a_y^{-1}\gamma(y)^2$, whence the idea of taking $a = \nu^{2}$.

\subsubsection{Conjugate gradient improvement and stabilization}

Adding a penalization also allows to accelerate the conjugate gradient algorithm, indeed it reduces the conditioning of the Hessian matrix by killing the small eigenvalues, and therefore accelerates the conjugate gradient algorithm's convergence. It also stabilizes this algorithm, indeed when $\eps$ is small we observe that without penalization, the numerical error may cause instabilities by returning a non positive definite Hessian. Adding the positive definite Hessian of the penalization function bypasses this instability.

\subsection{Epsilon scaling}\label{subsect:eps_scaling}

For all entropic algorithms, we observe that when $\eps$ is small, the algorithm may be very slow to find the region of optimality. For the Bregman projection, the formula for the speed of convergence in Subsection \ref{subsect:converge_Bregman} suggests to have a strategy of $\eps-$scaling: i.e. we solve the problem for $\eps = 1$, so that the function to optimize is very smooth. Then solve the problem for $\eps' < \eps$, with the previous optimum as an initial point. We continue this algorithm until we reach the desired value for $\eps$. In practice we divide $\eps$ by $2$ at each step.

\subsection{Grid size adaptation}

It may be a huge loss of time to run the algorithm on full resolution since the beginning of $\eps-$scaling. To prevent this waste of time, Schmitzer \cite{schmitzer2016stabilized} suggests to raise the size of the grid at the same time than shrinking $\eps$. In practice we give to each new point of the grid for $\varphi$, $\psi$, and $h$ the value of the closest point in the previous grid. We use heuristic criteria to decide when to doble the size of the grid, avoiding for example to doble is when $\eps$ is too small as is seriously challenges the stability of the resolution scheme.

\subsection{Kernel truncation}

While $\eps$ shrinks to $0$, we observe that the optimal transport tends to concentrate on graphs, as suggested in \cite{de2018local}. Because of the exponential, the value of the optimal probability far enough to these graphs tends to become completely negligible. For this reason, Schmitzer \cite{schmitzer2016stabilized} suggests to truncate the grid in order to do much less calculation. In dimensions higher than $1$, the gain in term of number of operation may quickly reach a factor $100$ for small $\eps$. In practice we removed the points in the grid when their probability were smaller than $10^{-7}\mu_x$ (resp. $10^{-7}\nu_y$) for all $x\in\Xc$ (resp. for all $y\in\Yc$).

\subsection{Computing the concave hull of functions}\label{subsect:convex_hull}

We were not able to find algorithms that compute the concave hull of a function in the literature, so we provide here the one we used. Let $f:\Yc\longrightarrow\R$.

In dimension $1$ the algorithm is linear in $|\Yc|$, we use the McCallum \& Avis \cite{mccallum1979linear} algorithm to find the points of the convex hull of the upper graph of $f$ in a linear time and then we go through these points until we find the two consecutive points $y_1,y_2\in\Yc$ around the convex hull such that $y_1 < x \le y_2$. Then $f_{conc}(x) = \frac{y_2-x}{y_2-y_1}f(y_1)+\frac{x-y_1}{y_2-y_1}f(y_2)$.

\begin{algorithm}
\caption{Concave hull of $f$.}\label{algo:convex}
\begin{algorithmic}[1]
\Procedure{ConcaveHull}{$f,x,grid,gradientGuess$}
\If {$gradientGuess$ is ${\rm None}$}
\State $grad \gets \text{vector of zeros with the same size than }x$
\State $gridF \gets f(grid)$
\Else
\State $grad \gets gradientGuess$
\State $gridF \gets f(grid)-grad\cdot grid$
\EndIf
\State $y \gets \argmax gridF$
\State $support \gets [y]$
\State $gridF \gets gridF - gridF[y_0]$
\While {${\rm True}$}
\If {$x\in\aff\, support$}
\State $bary\gets \text{barycentric coefficients of }x\text{ in the basis }support$
\If {bary\text{ are all }> 0}
\State $value \gets \text{sum}\big(bary\x f(support)\big)$
\State\Return $\{\text{"value"} : value ; \text{"support"}:support ;$\par
        \hskip\algorithmicindent\quad \quad\quad\quad $\text{"barycentric coefficients"}:bary ; \text{"gradient"}:grad\}$
\Else
\State $i\gets \argmin bary$
\State $\text{remove entry }i\text{ in }support$
\State $\text{remove entry }i\text{ in }bary$
\EndIf
\Else
\State $projx \gets \text{orthogonal projection of }x\text{ on }\aff\,support$
\State $p = x-projx$
\State $scalar \gets p\cdot(grid-x)$
\If {$scalar \text{ are all }\le 0$}
\State Fail with error "x not in the convex hull of grid."
\EndIf
\State $y\gets \argmax\{gridF/scalar\text{ such that }scalar >0\}$
\State $\text{add }y\text{ to }support$
\State $a\gets -gridF[y]/scalar[y]$
\State $gridF\gets gridF+a\x scalar$
\State $grad \gets grad - a\x p$
\EndIf
\EndWhile
\EndProcedure
\end{algorithmic}
\end{algorithm}

In higher dimension we use Algorithm \ref{algo:convex} in order to compute the convex hull of a function. We do not know if a better algorithm exists, but this one should be the fastest when the active points of the convex hull are already close to the maximum, this will be useful to compute $(c(x,\cdot) - \psi)_{conc}(x)$ from Theorem \ref{thm:entropy_error} below, so as the field "gradient" of the result that allows to find the right $h$. We believe that the complexity of this algorithm is quadradic in the (not so improbable) worst case of a concave function, $O\big(n\ln(n)\big)$ on average for a "random" function, and linear when the guess of the gradient is good. These assesments are formal and based on the observation of numerics, we do not prove anything about Algorithm \ref{algo:convex}, not even the fact that it cannot go on infinite loops. We provide it for the reader who would like to reproduce the numerical experiments without having to search for an algorithm by himself.

\section{Convergence rates}\label{sect:conv_rate}

\subsection{Discretization error}

\begin{Proposition}\label{prop:approx_marg}
Let $\mu\preceq\nu$ in convex order in $\Pc(\R^d)$ having a dual optimizer $(\varphi,\psi,h)\in\Dc_{\mu,\nu}(c)$ such that $\varphi$ is $L_\varphi-$Lipschitz, and $\psi$ is $L_\psi-$Lipschitz. Then for all $\mu'\preceq\nu'$ in convex order in $\Pc(\R^d)$ having a dual optimizer $(\varphi',\psi',h')\in\Dc_{\mu,\nu}(c)$ such that $\varphi'$ is $L_{\varphi'}-$Lipschitz, and $\psi'$ is $L_{\psi'}-$Lipschitz, we have
\b*
\left|\Sbf_{\mu,\nu}(c)-\Sbf_{\mu',\nu'}(c)\right|\le \max\left(L_\varphi,L_{\varphi'}\right) W_1(\mu',\mu)+\max\left(L_\psi,L_{\psi'}\right) W_1(\nu',\nu)
\e*
\end{Proposition}
The proof of Proposition \ref{prop:approx_marg} is reported in Subsection \ref{subsect:approx_marg}.

\begin{Remark}
Guo \& Ob{\l}{\'o}j \cite{guo2017computational} provide a very similar result in Proposition 2.2. It is however very different as they need to introduce an approximately martingale optimal transport problem, and our result makes hypotheses on the Lipschitz property of the dual optimizers, which are unknown, and even their existence in unknown. In dimension $1$, thanks to the work by Beiglb\"ock, Lim \& Ob{\l}{\'o}j \cite{beiglbock2017dual}, we may prove the existence of these Lipschitz dual, thanks to some regularity assumptions on $c$. In higher dimension, there are ongoing investigations about the existence of similar results. However, by Example 4.1 in \cite{de2018quasi}, it will be necessary to make assumptions on $\mu,\nu$ as well, as the smoothness of $c$ cannot guarantee the existence of a dual optimizer. Proposition \ref{prop:approx_marg} is entitled to be a proposition of practical use, we may formally assume that the partial dual functions that we get converge to the continuous dual and assume that their Lipschitz constant converges to the Lipschitz constant of the limit.
\end{Remark}

We refer to Subsection 2.2 in \cite{guo2017computational} for a study of the discrete $\Wc_1-$approximation of the continuous marginals. In dimensions higher than $3$, it is necessary to use a Monte-Carlo type approximation of $\mu$ and $\nu$ to avoid the curse of dimensionality linked to a grid type approximation. However, Proposition \ref{prop:approx_marg} is not well-adapted to estimate the error, as we know from \cite{fournier2015rate} that the Wasserstein distance between a measure and its Monte-Carlo estimate is of order $n^{\frac{-1}{d}}$. Next proposition deals with this issue. For two sequences $(u_N)_{N\ge 0}$ and $(v_N)_{N\ge 0}$, we denote $u_N\approx v_N$ when $N\longrightarrow \infty$ if $u_N/v_N$ converges to $1$ in probability, when $N\longrightarrow \infty$.

\begin{Proposition}\label{prop:approx_marg_MC}
Let $\mu\preceq\nu$ in convex order in $\Pc(\R^d)$ having a dual optimizer $(\varphi,\psi,h)\in\Dc_{\mu,\nu}(c)$, and $\mu_N$ and $\nu_M$, independent Monte-Carlo estimates of $\mu$ and $\nu$ with $N$ and $M$ samples. If furthermore $\mu'_N\preceq\nu'_M$ is in convex order in $\Pc(\R^d)$ having a dual optimizer $(\varphi_N,\psi_M,h')\in\Dc_{\mu'_N,\nu'_M}(c)$ such that when $N,M\longrightarrow\infty$, we have

\no{\rm (i)} $(\mu_N-\mu)[\varphi]\approx (\mu_N-\mu)[\varphi_N]$,
\quad
\no{\rm (ii)} $(\nu_M-\nu)[\psi]\approx (\nu_M-\nu)[\psi_M]$,

\no{\rm (iii)} $(\mu_N-\mu'_N)[\varphi]\approx (\mu_N-\mu'_N)[\varphi_N]$,
\quad
\no{\rm (iv)} $(\nu_M-\nu'_M)[\psi]\approx (\nu_M-\nu'_M)[\psi_M]$.

Then we have
\b*
\left|\Sbf_{\mu,\nu}(c)-\Sbf_{\mu',\nu'}(c)\right|\le \a \sqrt{\frac{{\rm Var}_\mu[\varphi]}{N}+\frac{{\rm Var}_\nu[\psi]}{M}}+\big|(\mu_N-\mu'_N)[\varphi_N]\big|+\big|(\nu_M-\nu'_M)[\psi_M]\big|,
\e*
with probability converging to $1-2\int_\a^\infty \exp\left(-x^2/2\right)dx$, when $N,M\longrightarrow\infty$.
\end{Proposition}
The proof of Proposition \ref{prop:approx_marg_MC} is reported in Subsection \ref{subsect:approx_marg}.

\begin{Remark}
In Proposition \ref{prop:approx_marg_MC}, we introduce $\mu_N',\nu_M'$ because the Monte-Carlo approximation will not conserve the convex order for $\mu_N,\nu_N$ in general. Then we obtain $\mu_N',\nu_M'$ from "convex ordering processes" such as the one suggested in Subsection \ref{subsubsect:convex_order}, or the one suggested in \cite{alfonsi2017sampling}. In both cases, the quantity $\big|(\mu_N-\mu'_N)[\varphi_N]\big|+\big|(\nu_M-\nu'_M)[\psi_M]\big|$ may be computed exactly numerically.
\end{Remark}

\subsection{Entropy error}

In this subsection, $m_\eps$ is a generic finite measure and no assumptions are made on $\mu_\eps$ and $\nu_\eps$.

\begin{Theorem}\label{thm:entropy_error}
Let $\mu_\eps\preceq\nu_\eps$ in convex order in $\Pc(\R^d)$ with dual optimizer $(\varphi_\eps,\psi_\eps,h_\eps)\in\Dc_{\mu_\eps,\nu_\eps}(c)$ to the $\eps-$entropic dual problem with reference measure $m_\eps$, such that we may find $\gamma,\eta,\beta>0$, sets $(D_\eps^\Xc,D_\eps^\Yc)_{\eps>0}\subset \Bc(\R^d)$, and parameters $(\a_\eps,A_\eps)_{\eps>0}\subset(0,\infty)$, such that if we denote $r_\eps:= \eps^{\frac12-\eta}$, $m_\eps^\Xc:=m_\eps\circ X^{-1}$, and $\Delta_\eps:= \varphi\oplus\psi+h^\otimes-c$, for $\eps >0$ small enough we have:

\no{\rm (i)} $\frac{d\mu_\eps}{dm^{\Xc}_\eps}(x)\le \eps^{-\gamma}$, $\mu_\eps-$a.e., $m_\eps[\Omega]\le \eps^{-\gamma}$, and $A_\eps\ll \eps^{-q}$ for all $q>0$.

\no{\rm (ii)} For $m_\eps-$a.e. $(x,y)\in \R^d\x\R^d$, we have $(m_{\eps})_{x}\left[B_{\a_\eps}(y)\right]\ge \eps^\gamma$, and for $(m_\eps)_x-$a.e. $y'\in B_{\a_\eps}(y)$, we have $|\Delta_\eps(x,y)-\Delta_\eps(x,y')|\le \gamma\eps\ln(\eps^{-1})$.

\no{\rm (iii)} $\mu_\eps\left[(D_\eps^{\Xc})^c\right]\ll 1/\ln(\eps^{-1})$ and for all $x\in D^{\Xc}_\eps$ we may find $k_x^\eps\in \N$, $S_x^\eps\in (B_{A_\eps})^{k_x^\eps}$, and $\lambda_x^\eps\in [0,1]^{k_x^\eps}$ with $\det_\aff(S_x^\eps)\ge A_\eps^{-1}$, $\min\lambda_x^\eps\ge A_\eps^{-1}$, and $\sum_{i=1}^{k_x^\eps}\lambda_{x,i}^\eps S_{x,i}^\eps= x$, convex combination.

\no{\rm (iv)} On $B_{r_\eps}(S_x^\eps)$, $\Delta_\eps(x,\cdot)$ is $\Ctn^2$, $A_\eps^{-1} I_d\le\partial^2_y\Delta_\eps(x,\cdot) \le A_\eps I_d$, and for all $y,y'\in B_{r_\eps}(S_x^\eps)$, we have $|\partial^2_y\Delta_\eps(x,y)-\partial^2_y\Delta_\eps(x,y')|\le \eps^\eta$.

\no{\rm (v)} For $x\in D_\eps^\Xc$ and $y\notin B_{r_\eps}(S_x^\eps)$, we have that $\Delta_\eps(x,y)\ge \sqrt{\eps}\, \dist(y,S_x^\eps)$.

\no{\rm (vi)} $\nu_\eps\left[(D_\eps^\Yc)^c\right]\ll \frac{1}{A_\eps\ln(\eps^{-1})}$ and for all $y_0\in D_\eps^\Yc$, $R,L \ge 1$, and $f:\R^d\longrightarrow \R_+$ such that $\frac{f}{\|f\|_\infty^R}$ is $L-$Lipschitz, we have
\b*
\left|\int_{B_{R}}f(y)\left[\frac{d(m_\eps)_{x}\circ\zoom_{\sqrt{\eps}}^{y_0}}{(m_{\eps})_{x}[B_{R\sqrt{\eps}}(y_0)]}-\frac{dy}{|B_{R}|}\right]\right|\le [R+L]^\gamma\eps^{\beta}\int_{B_{R}}f(y)\frac{dy}{|B_{R}|}.
\e*
Then if we denote $\P_\eps := e^{-\frac{\Delta_\eps}{\eps}}m_\eps$, we have
\b*
\lim_{\eps\to 0}\frac{\mu_\eps[\bar\varphi_\eps]+\nu_\eps[\psi_\eps]-\P_\eps[c]}{\eps} =\frac{d}{2},&\mbox{where}&\bar\varphi_\eps := \big(c(X,\cdot)-\psi_\eps\big)_{conc}(X).
\e*
\end{Theorem}

The proof of Theorem \ref{thm:entropy_error} is reported to Subsection \ref{subsect:entropy_error}.

\begin{Corollary}
Under the assumptions of \ref{thm:entropy_error}, we have that $\P_\eps[c]\ge \Sbf_{\mu,\nu}(c)-\frac{d}{2}\eps + o(\eps)$, when $\eps\longrightarrow 0$.
\end{Corollary}
\proof
We fix $h(x)\in \partial (c(x,\cdot)-\psi_\eps)_{conc}(x)$ for all $x\in\Xc$, then
$$\big((c(X,\cdot)-\psi_\eps)_{conc}(X),\psi_\eps,h\big)\in \Dc_{\mu_\eps,\nu_\eps}(c),$$
and therefore $\mu_\eps\big[(c(X,\cdot)-\psi_\eps)_{conc}(X)\big]+\nu_\eps[\psi_\eps]\ge \Ibf_{\mu_\eps,\nu_\eps}(c) \ge \Sbf_{\mu_\eps,\nu_\eps}(c)\ge \P_\eps[c]$. Theorem \ref{thm:entropy_error} concludes the proof.
\ep\\

\begin{figure}[h]
\centering
\begin{subfigure}{.49\textwidth}
  \centering
  \includegraphics[width=.99\linewidth]{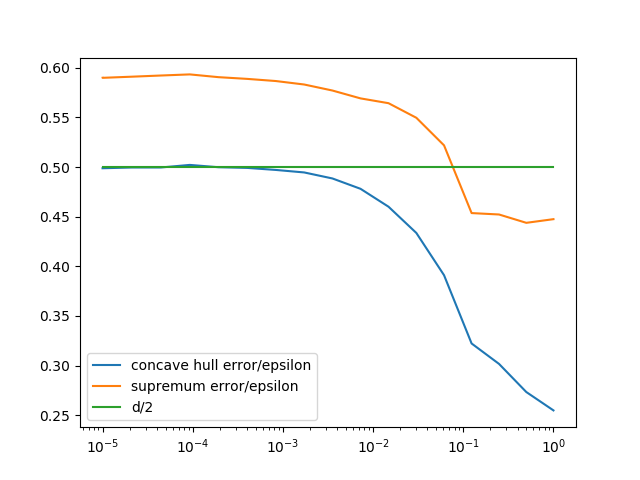}
  \caption{Dimension 1.}
  \label{fig:errordim1}
\end{subfigure}%
\begin{subfigure}{0.49\textwidth}
  \centering
  \includegraphics[width=.99\linewidth]{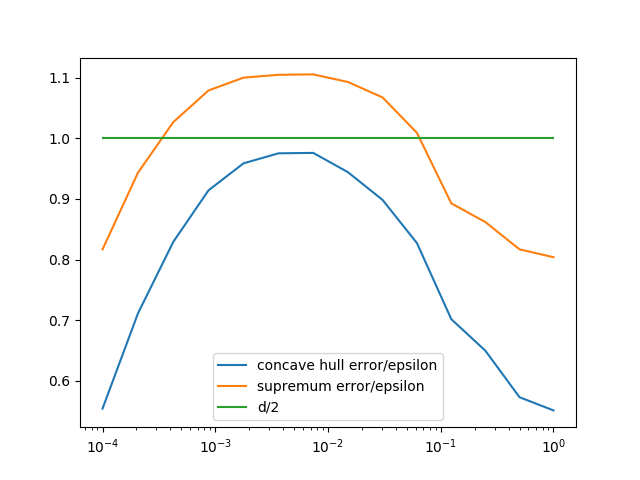}
  \caption{Dimension 2.}
  \label{fig:errordim2}
\end{subfigure}
\caption{Duality gap for the supremum, and the concave hull dual approximation vs $\eps$.}
\label{fig:table of gaps}
\end{figure}

Figure \ref{fig:table of gaps} gives numerical examples of the convergence of the duality gap when $\eps$ converges to $0$. In these graphs, the blue curve gives the ratio of the dominator of the duality gap $\mu_\eps\Big[\big(c(X,\cdot)-\psi_\eps\big)_{conc}(X)\Big]+\nu_\eps[\psi_\eps]-\P_\eps[c]$ with respect to $\eps$. It is meant to be compared to the green flat curve which is its theoretical limit $\frac{d}{2}$ according to Theorem \ref{thm:entropy_error}. Finally the orange curve provides the ratio of the weaker dominator $\mu_\eps\Big[\sup\big(c(X,\cdot)-\psi_\eps\big)\Big]+\nu_\eps[\psi_\eps]-\P_\eps[c]$ of the duality gap with respect to $\eps$. The interest of this last weaker dominator is that it avoids computing the concave envelop which may be a complicated issue, while having a reasonable comparable performance in practice than the concave hull dominator as showed by the graphs and by Remark \ref{rmk:onlysup} below.

Figure \ref{fig:errordim1} provides these curves for the one-dimensional cost function $c := XY^2$, $\mu$ uniform on $[-1,1]$, and $\nu = |Y|^{1.5}\mu$. The grid size adaptation method is used and the size of the grid goes from $10$ when $\eps = 1$ to $10000$ when $\eps = 10^{-5}$. Figure \ref{fig:errordim2} provides these curves for the two-dimensional cost function $c:(x,y)\in\R^2\x\R^2\longmapsto x_1(y_1^2+2y_2^2)+x_2(2y_1^2+y_2^2)$, $\mu$ uniform on $[-1,1]^2$, and $\nu = (|Y_1|^{1.5}+|Y_2|^{1.5})\mu$. The grid size adaptation method is used again and the size of the grid goes from $10\x 10$ when $\eps = 1$ to $160\x 160$ when $\eps = 10^{-4}$.

\begin{Remark}
The hypotheses of regularity are impossible to check at the current state of the art. One element that could argue in this direction is the fact that the limit of $\Delta_\eps$ when $\eps\longrightarrow 0$, that we shall denote $\Delta_0$, should satisfy the Monge-Amp\`ere equation $\det\left(\partial^2_y\Delta_0\right) = f$, similar to classical transport \cite{trudinger2006second} that may provide some regularity, but less than the one needed to satisfy (ii) of Theorem \ref{thm:entropy_error}, see Chapter 5 of \cite{gutierrez2001monge}. It also justifies, in the case where $f>0$, that $\partial_y^2\Delta_\eps\in GL_d(\R)$.
\end{Remark}

\begin{Remark}
Assumption (vi) on the local convergence of the reference measures is justified is we take a regular grid for $\Yc$ that becomes fine fast enough. If the grid does not become fine fast enough, then the local decrease of $\Delta_\eps-\Delta_\eps(x)-\nabla \Delta_\eps(x)\cdot(Y-x)$ is exponential because of the shape of the kernel. We observe on the experiments that in this case the duality gap becomes indeed smaller, however as a downside, the convergence of the scheme becomes much less effective because the marginal error stays high (see Proposition \ref{prop:approx_marg}). See Figure \ref{fig:errordim2}.
\end{Remark}

\begin{Remark}
Assumption (ii) from the theorem seems to hold in one dimension, but it seems to be wrong in two dimensions, as shown in the example of Figure \ref{fig:test}. However, we may still find a formula similar to \eqref{eq:epsilon_times} below. Therefore, we may reasonably assume that the error is still of order $\eps$, as confirmed in a the numerical examples by Figure \ref{fig:errordim2}.
\end{Remark}

\begin{Remark}
In the case when $\Xc$ is obtained from Monte-Carlo methods, (vi) is verified with a constant that depends on the point, and is probabilistic. the exponent $\alpha$ may be taken taken equal to $\frac{1}{d}$, see \cite{fournier2015rate}.
\end{Remark}

\begin{Remark}
Despite the difficulty to check the assumptions, this result is inspired and satisfied by observation on the numerics, we have tried with several cost functions, differentiable or not, and the result seems to be always satisfied, probably with the help of its universality. See figure \ref{fig:table of gaps}.
\end{Remark}

\begin{Remark}
An easier version of Theorem \ref{thm:entropy_error} could be obtained by the same strategy for classical optimal transport under the twist condition $c_x(x,\cdot)$ injective, replacing $\big(c(X,\cdot)-\psi_\eps\big)_{conc}(X)$ by $\sup\big(c(X,\cdot)-\psi_\eps\big)$.
\end{Remark}

\begin{Remark}\label{rmk:onlysup}
By similar reasoning, we may prove that even for martingale optimal transport, replacing $\big(c(X,\cdot)-\psi_\eps\big)_{conc}(X)$ by $\sup\big(c(X,\cdot)-\psi_\eps\big)$ would still give a good result (see \ref{fig:table of gaps}). We may formally estimate the new limit of $\frac{duality\,gap}{\eps}$: for all $y_i^\eps(x)$ that are not the optimizer (say $y_0^\eps(x)$), the weight added to $\frac{d}{2}$ is $\lambda_i^\eps(x)\big(\Delta_\eps(x,y_i^\eps(x))-\Delta_\eps(x,y_0^\eps(x))\big)$. By using the tools of the proof of Theorem \ref{thm:entropy_error}, we get the formal formula, if we denote $\lambda_i$ for the limit of $\lambda_i^\eps(X)$, and $y_i$ for the limit of $y_i^\eps(x)$, we have
\be\label{eq:epsilon_times}
\mu_\eps\big[\sup \{c(X,\cdot)-\psi_\eps\}(X)\big]+\nu_\eps[\psi_\eps]-\P_\eps[c]\approx \a\eps,
\ee
with $\a=\frac{d}{2}+\int_{\R^d}\sum_{i>0}\lambda_i\ln\left(\frac{dm^{\Yc}\big(y_i\big)}{dm^{\Yc}\big(y_0\big)}\frac{\lambda_0\det\left(\partial^2_y\Delta_0(x,y_0)\right)}{\lambda_i\det\left(\partial^2_y\Delta_0(x,y_i)\right)}\right)d\mu$.
Then we could reasonably make the assumption the the second term in $\a$ does not explode, and then the limit is still of the order of $\eps$, as we may see on the numerical experiments of Figure \ref{fig:table of gaps}. This result also generalises to optimal transport when there are several transport maps.
\end{Remark}

\subsection{Penalization error}

\begin{Proposition}\label{prop:limit_pen}
Let $(\mu,\nu)\in \Pc(\Xc)\x\Pc(\Yc)$ in convex order. Let $\nu_\alpha:=\P_\alpha\circ Y^{-1}$, where $\P_\alpha$ is the optimal probability for the entropic dual implied problem with an additional penalization $\alpha f$, where $f$ is a super-linear, strictly convex, and differentiable function $\R^\Yc\longrightarrow\R$. Then let $\psi_{0}$ be the only optimizer for the entropic dual implied problem with minimal $f(\psi)$, we have
\b*
\frac{\nu_\alpha-\nu}{\alpha}\longrightarrow \nabla f(\psi_{0}),&\mbox{when}&\alpha\longrightarrow 0.
\e*
\end{Proposition}

The proof of Proposition \ref{prop:limit_pen} is reported to Subsection \ref{subsect:limit_pen}.

\subsection{Convergence rates of the algorithms}

\subsubsection{Convergence rate for the simplex algorithm}\label{subsubsect:converg_simplex}

Precise results on the convergence rate of the simplex algorithm is an open problem. Roos \cite{roos1990exponential} gave an example in which the convergence takes an exponential time in the number of parameters. However the simplex algorithm is much more efficient in practice, Smale \cite{smale1983average} proved that in average, the number of necessary steps in polynomial in the number of entries, and Spielman \& Teng \cite{spielman2004smoothed} refined this analysis by including the number of constraints in the polynomial.

However, none of these papers provide the real time of convergence of this algorithm. Schmitzer \cite{schmitzer2016stabilized} reports that this algorithm is not very useful in practice as it only allows to solve a discretized problem with no more than hundreds of points.

\subsubsection{Convergence rate for the semi-dual algorithm}

We notice that any subgradient of this function is a difference of probabilities, and then the gradient is bounded. Furthermore the function $V$ is a supremum of a finite number of affine functions, and therefore it does not have a smooth second derivative. In this condition the best theoretical way to optimize this function is by a gradient descent with a step size of order $O(1/\sqrt{n})$ at the $n-$th step, see Ben-Tal \& Nemirovski \cite{ben2001lectures}. Then by Theorem 5.3.1 of \cite{ben2001lectures}, the rate of convergence is $O(1/\sqrt{n})$ as well, which is quite slow. Furthermore, the time of computation of one step needs to compute one convex hull which has the average complexity $O(|\Yc|\ln(|\Yc|))$ for each $x\in\Xc$, and $O(|\Yc|)$ in dimension 1, see Subsection \ref{subsect:convex_hull}. However, we give in Subsection \ref{subsect:convex_hull}) an algorithm that computes the concave hull in a linear time if the relying points of the concave hull do not change too much. Then let us be optimistic and assume that the computation of one concave hull is on average $O(|\Yc|)$, then we have that the complexity is $O(|\Xc||\Yc|)$ operations for each step. Although this algorithm is highly parallelizable, its complexity imposes to the grid to be very coarse. Indeed, to get a precision of $10^{-2}$, we need an order of $10^4$ operations. We shall see that the entropic algorithms are much more performing for this low precision.

Notice that even though the best theoretical algorithm is the last gradient descent, Lewis \& Overton \cite{lewis2013nonsmooth} showed that in most case, quasi-Newton methods converge faster, even when the convex function is non-smooth, however they find a particular case in which quasi-Newton fails at being better. The L-BFGS method is a quasi-Newton method that is adapted to high-dimensions problems. The Hessian (even though it does not exist) is approximated by a low-dimensional estimate, and the classical Newton step method is applied. See \cite{wright1999numerical} for the exact algorithm. We see on simulations that this algorithm is indeed much more efficient.


Even if the quasi-Newton algorithm gives better results, the smooth entropic algorithms are much more effective in practice.

\subsubsection{Convergence rate for the Sinkhorn algorithm}\label{subsect:converge_Bregman}

In practice, if we want to observe the transport maps from \cite{de2018local} to have a good precision on the estimation of the support of the optimal transport, we need to set $\eps = 10^{-4}$. The rate of convergence of the Sinkhorn algorithm is given by $\kappa^{2n}$ after $n$ iterations for some $0<\kappa<1$, see \cite{knight2008sinkhorn}. This result is extended by \cite{guo2017computational} to the one-dimensional martingale Sinkhorn algorithm. However, we have $\kappa:=\frac{\sqrt{\theta}-1}{\sqrt{\theta}+1}$, and
$$\theta:=\max_{(x_1,y_1),(x_2,y_2)\in\Omega}\exp\left(\frac{c(x_1,y_1)+c(x_2,y_2)-c(x_2,y_1)-c(x_1,y_2)}{\eps}\right),$$
in the case of classical transport. Then $\theta$ is of the order of $\exp(K(c)/\eps)$ for some map $K(c)$ bounded from below. For $\eps = 10^{-5}$, this $\theta$ is so big that $\kappa^2$ is so close to $1$ that $\kappa^{2n}$, with $n$ the number of iterations will remain approximately equal to $1$.

We also see in practice for the martingale Sinkhorn algorithm that the rate of convergence is not exponential for small values of $\eps$, see Figure \ref{fig:performances} in the numerical experiment part, as the graph is logarithmic in the error, an exponential convergence rate would be characterized by a straight line. However we observe that for the Bregman projection algorithm we do not have a straight line during the first part of the iteration for the one-dimensional case, and it never happens in the two-dimensional case.

In this regime of $\eps$ small, another convergence theory looks to have a better fit with this algorithm. The Sinkhorn algorithm may be interpreted as a block coordinates descent for the optimization of the map $V_\eps(\varphi,\psi)$. We optimize alternatively in $\varphi$, and in $\psi$. We know from Beck \& Tetruashvili \cite{beck2013convergence} that this optimization problem has a speed of convergence given by $\frac{LR(x_0)^2}{n}$, where $R(x_0)$ is a quantity that is of the order of $|x_0-x^*|$ in practice, where $x^*$ is the closest optimizer of $V_\eps$ and $L$ is the Lipschitz constant of the gradient. This speed is more comparable to the convergence observed in practice. More precisely, $L$ is of the order of $1/\eps$. This formula shows that in order to minimize the problem for a very small $\eps$, we first need to make $R(x_0)$ small to compensate $L$. This can be done by minimizing the problem for larger $\eps$. In practice, we divide $\eps$ by $2$ until we reach a sufficiently small $\eps$. Then we make the grid finer as $\eps$ becomes small, and exploit the sparsity in the problem that appears when $\eps$ gets small. See Schmitzer \cite{schmitzer2016stabilized}.

We may apply the same theory for the martingale $V_\eps$ and its block optimization in $(\varphi,h)$ and in $\psi$. Let $\Dc_{\Xc,\Yc} := \{(\varphi,\psi,h)\in \R^\Xc\x\R^\Yc\x(\R^d)^\Xc\} \approx \R^{(d+1)|\Xc|+|\Yc|}$, and for $x := (\varphi,\psi,h)\in\Dc_{\Xc,\Yc}$, let $\Delta(x) := (\varphi\oplus\psi+h^\otimes)_{\Xc\x\Yc}$.

\begin{Theorem}\label{thm:cvg_Sinkhorn}
Let $x_0 = (\varphi_0,\psi_0,h_0)\in \Dc_{\Xc,\Yc}$ such that $\left(e^{-\frac{\varphi_0\oplus\psi_0+h_0^\otimes-c}{\eps}}\right)_{\Xc\x\Yc}$ sums to $1$, and for $n\ge 0$, let the $n^{th}$ iteration of the martingale Sinkhorn algorithm:
\b* 
x_{n+1/2} &:=& \left(\varphi_n,\psi_{n+1}:=\argmin_\psi V_\eps(\phi_n,\cdot,h_n),h_n\right),\\
x_{n+1} &:=& \left(\varphi_{n+1}:=\argmin_\varphi V_\eps(\cdot,\psi_{n+1},\cdot),\psi_{n+1},h_{n+1}:=\argmin_h V_\eps(\cdot,\psi_{n+1},\cdot)\right).
\e*
Furthermore let $\P_0\in\Mc(\mu,\nu)$ and let $\Xc^*$ be the minimizing affine space of $V_\eps$ and let $V_\eps^*$ be its minimum, then we have
\be
V_\eps(x_n)-V_\eps^*&\le&
\frac{\beta R(x_0)^2\eps^{-1}}{n},\label{eq:cvg_linear}\\
V_\eps(x_n)-V_\eps^*&\le&
\left(1-\frac{2|\Xc| }{\beta\lambda_{2}}e^{-\frac{D(x_0)}{\eps}} \right)^{n}\big(V_\eps(x_0)-V_\eps^*\big),\label{eq:cvg_harm}\\
\mbox{and}\quad\dist(x_n,\Xc^*)&\le& \sqrt{\frac{\lambda_{2}\eps}{|\Xc|} e^{\frac{D(x_0)}{\eps}}} \big(V_\eps(x_n)-V_\eps^*\big)^{\frac12},\nonumber
\ee
where $\beta := \big(2\ln(2)-1\big)^{-1}\approx 2,6$,\quad $\lambda_p :=
|\Xc|^{-1}\underset{x\in\Dc_{\Xc,\Yc}:\Delta(x)\neq 0}{\inf}\,\underset{\Delta(\widetilde{x})=\Delta(x)}{\sup}\frac{|\Delta(x)|_p^p}{|\widetilde{x}|_p^p}$, for $p = 1,2$,

$D(x_0):= \lambda_1\max(1,\|\Yc-\Xc\|_\infty)\frac{V_\eps(x_0)-\P_0[c]-\eps}{|\Xc|(\P_0)_{\min}}$,


$(\P_0)_{\min}:=\min_{x\in\Xc\x\Yc}\P_0[\{x\}]$,\quad $\|\Yc-\Xc\|_\infty:=\sup_{x\in\Yc-\Xc}|y-x|_\infty$,
\b*
\mbox{and we may choose $R(x_0)$ among}& R(x_0):=\left\lbrace
\begin{array}{c}
\sup_{V_\eps(x)\le V_\eps(x_0)}\dist(x,\Xc^*)\\
\sqrt{\frac{\lambda_{2}\eps}{|\Xc|} e^{\frac{D(x_0)}{\eps}}} \big(V_\eps(x_0)-V_\eps^*\big)^{\frac12}\\
2 \lambda_1\frac{V_\eps(x_0)-\P_0[c]-\eps}{|\Xc|(\P_0)_{\min}}\\
\sup_{k\ge 0}\dist(x_k,\Xc^*)
\end{array}\right. .
\e*
\end{Theorem}

The proof of Theorem \ref{thm:cvg_Sinkhorn} is reported in Subsection \ref{subsect:cvg_Sinkhorn}.

\begin{Remark}
By the same arguments, we may prove a similar theorem for the case of optimal transport with $\lambda_{1} = 
\min\left(1,\frac{|\Yc|}{|\Xc|}\right)$, and $\lambda_2 = \frac{1+\left(\frac{|\Yc|}{|\Xc|}\right)^2}{1+\frac{|\Yc|}{|\Xc|}}$.
\end{Remark}

\begin{Remark}
The theoretical rate of convergence given by Theorem \ref{thm:cvg_Sinkhorn} becomes pretty bad when $\eps\longrightarrow 0$. We observe it in practise when we apply this algorithm with a small $\eps$ and a starting point $x_0 = 0$. This emphasizes the need of using the epsilon scaling trick, of Subsection \ref{subsect:eps_scaling}.
\end{Remark}

\begin{Remark}
Depending on the experiment, in some cases we observe a linear convergence like \eqref{eq:cvg_linear} (see Figure \ref{fig:performancesdim1}), however in other cases, we observe a convergence speed that looks more like \eqref{eq:cvg_harm} (see Figure \ref{fig:performancesdim2}). However, the convergence rates that we provide here are generic, if we wanted to have convergence rates that look more like the one observed, we would need to look for the asymptotic convergence rates like it was suggested by Peyr\'e in \cite{peyre2017computational} for the case of classical transport.
\end{Remark}

\begin{Remark}
The positive probability $\P_0\in \Mc(\mu,\nu)$ is necessary. We know from \cite{de2017irreducible} that for some (possibly elementary) $\mu \preceq \nu$ in convex order, we way find $(x_0,y_0)\in \Xc\x\Yc$ such that $\P[\{(x_0,y_0)\}] = 0$ for all $\P\in\Mc(\mu,\nu)$ even thought $\mu[\{x_0\}]>0$ and $\nu[\{y_0\}]>0$ (see Example 2.2 in \cite{de2017irreducible}). Therefore, in this situation there is no optimal $x^*\in \Dc_{\Xc\x\Yc}$, as this would mean that $\Delta(x^*)_{x_0,y_0} = \infty$.
\end{Remark}

\subsubsection{Convergence rate for the Newton algorithm}

When the current point gets close enough from the optimum, the convergence rate of the Newton algorithm is quadratic if the hessian is Lipschitz, i.e. $|x_k-x^*|$ and $|\nabla V_\eps(x_k)|$ both converge quadratically to $0$, see Theorem 3.5 in \cite{wright1999numerical}. The truncated Newton is a bit slower, but still has a superlinear convergence rate, see Theorem 7.2 in \cite{wright1999numerical}.

\subsubsection{Convergence rate for the implied Newton algorithm}

The important parameter for Newton algorithm is the Lipschitz constant of the Hessian of the objective function. However in the case of variable implicitation, the presence of $\partial_y^2 F^{-1}$ in the Hessian, and the addition of the variation of $y(x)$ in the Lipschitz analysis may kill the Lipschitz property of the Hessian of $\tilde{F}$. The following proposition solves this problem.
\begin{Proposition}\label{prop:equivNewton}
Let $(\tilde{x}_n)_{n\ge 0}$ the Newton iterations applied to $\tilde{F}$ starting from $\tilde{x}_0 := x_0\in\Xc$. Now let $(x_n)_{n\ge 0}$ the sequence defined by recurrence by $x_0 := x_0$, then for all $n\ge 0$, $y_n:= y(x_n)$, and let $(x,y)$ be the result of a Newton step from $(x_n,y_n)$, and we set $x_{n+1} := x$.

Then $(\tilde{x}_n)_{n\ge 0} = (x_n)_{n\ge 0}$.
\end{Proposition}
The proof of Proposition \ref{prop:equivNewton} is reported to Subsection \ref{subsect:equivNewton}. This proposition implies that the theoretical convergence of the Newton algorithm on $F$ can be extended to the Newton algorithm applied to $\tilde{F}$, indeed the partial minimization in $y$ only decreases the distance from the current point to the minimum around the minimum. In practice we observe that the convergence for this implied algorithm is much faster and much more stable than the non-implied Newton algorithm.

\section{Proofs of the results}\label{sect:proofs}

\subsection{Minimized convex function}\label{subsect:conv_min}

\no {\bf Proof of Proposition \ref{prop:varimpl}} Let $x_1,x_2\in\Ac$, $y_1,y_2\in\Bc$, and $0\le \lambda\le 1$. We have
\b*
\lambda F(x_1,y_1)+(1-\lambda)F(x_2,y_2)&\ge& F\big(\lambda(x_1,y_1)+(1-\lambda)(x_2,y_2)\big)\\
&&+\alpha\frac{\lambda(1-\lambda)}{2}|(x_1,y_1)-(x_2-y_2)|^2\\
&\ge&\tilde{F}(\lambda x_1+(1-\lambda)x_2)+\alpha\frac{\lambda(1-\lambda)}{2}|x_1-x_2|^2.
\e*
By minimizing over $y_1$ and $y_2$, we get
$$\lambda \tilde{F}(x_1,y_1)+(1-\lambda)\tilde{F}(x_2,y_2)\ge \tilde{F}(\lambda x_1+(1-\lambda)x_2)+\alpha\frac{\lambda(1-\lambda)}{2}|x_1-x_2|^2,$$
which establishes the $\alpha-$convexity of $\tilde{F}$.

Now if we further assume that $\alpha >0$ and $F$ is $\Ctn^2$, $y\longmapsto F(x,y)$ is $\alpha-$convex, and therefore strictly convex and super-linear. Hence, there is a unique minimizer $y(x)$. Using the first order derivative condition of this optimum, we have $\partial_y F\big(x,y(x)\big) = 0$. By the fact that $\partial_y^2F$ is positive definite (bigger than $\alpha Id$ by $\alpha-$convexity), we may apply the local inversion theorem, which proves that $y(x)$ is $\Ctn^1$ in the neighborhood of $x$. We also obtain $\partial^2_{yx}F\big(x,y(x)\big)+\partial^2_y F\big(x,y(x)\big)\nabla y(x) = 0$, which gives the following expression of $\nabla y$:
$$\nabla y(x) = -\partial^2_{y}F^{-1}\partial^2_{yx}F\big(x,y(x)\big).$$

Now we may compute the derivatives of $\tilde{F}$. By definition, we have $\tilde{F}(x) = F\big(x,y(x)\big)$, then just differentiating this expression, we get
\b*
\nabla \tilde{F}(x) &=& \partial_x F\big(x,y(x)\big)+\partial_y F\big(x,y(x)\big)\nabla y(x)\\
&=&\partial_x F\big(x,y(x)\big),
\e*
where the second equality comes from the fact that $\partial_y F\big(x,y(x)\big) = 0$ because $y(x)$ is a minimizer. Finally we get the Hessian by deriving again this expression and injecting the value of $\nabla y(x)$.
\ep

\subsection{Limit marginal}\label{subsect:marginal_limit}

\no{\bf Proof of Theorem \ref{thm:marginal_limit}} Let $\alpha>0$, we are considering the following minimization problem:
\be\label{eq:pen_opt_duality}
\inf_\psi \widetilde{V}_\eps(\psi) + \alpha f(\psi)&=& \inf_{\varphi,\psi,h}\mu[\varphi]+\nu[\psi]+\eps\int e^{\frac{-\Delta}{\eps}}dm_0+\alpha f(\psi)\nonumber\\
&=&\inf_\psi\sup_{\P\in\Mc(\mu)}\P[c-\psi]-\eps H(\P|m_0)+\nu[\psi]+\a f(\psi)\\
&=&\sup_{\P\in\Mc(\mu)}\inf_\psi \P[c]-\eps H(\P|m_0) + (\nu-\P\circ Y^{-1})[\psi]+\a f(\psi)\nonumber\\
&=& \sup_{\P\in\Mc(\mu)}\P[c]-\eps H(\P|m_0)-(\a f)^*(\P\circ Y^{-1}-\nu)\nonumber\\
&=& \sup_{\P\in\Mc(\mu)}\P[c]-\eps H(\P|m_0)-\alpha^{-\frac{1}{p-1}}f^*(\P\circ Y^{-1}-\nu),
\ee
where the first equality comes from a mutualisation of the infima, the second comes from a partial dualisation of the infimum in $\varphi,h$ in a supremum over $\P\in\Mc(\mu)$, we obtain the third equality by applying the minimax theorem and reordering the terms, the fourth equality the definition of the Fenchel-Legendre transform, and the fifth and final equality is just a consequence of the transformation of a multiplyer of a $p-$homogeneous function by the Fenchel-Legendre conjugate. Let $(\alpha_n)_{n\ge 1}$ converging to $0$. As $\Yc$ is finite, the set $\Pc(\Yc)$ is compact. Then we may assume up to extracting a subsequence that $\nu_{\alpha_n}$ converges to some limit $\nu_l$. The first order optimality equation for all $y\in\Yc$ gives that $\nu - \nu_{\alpha_n} + \alpha_n \nabla f(\psi_n)$, where $\psi_n$ is the unique optimizer of $\widetilde{V}_\eps+\alpha f$. By the $p-$homogeneity of $f$, the gradient $\nabla f$ is $(p-1)-$homogeneous. Then we have the convergence $\tilde\psi_n :=\frac{\psi_{n}}{\alpha_n^{\frac{1}{p-1}}}\underset{n \to \infty}{\longrightarrow} \psi_l := \nabla f^{-1}(\nu_l-\nu)$. As we have by \eqref{eq:pen_opt_duality}
\b*
\inf_\psi \widetilde{V}_\eps(\psi) + \alpha_n f(\psi)=\inf_\psi\sup_{\P\in\Mc(\mu)}\P[c-\psi]-\eps H(\P|m_0)+\nu[\psi]+\a_n f(\psi)
\e*
By dividing this equation by $\alpha_n^{\frac{1}{p-1}}$, we have that $\psi_l$ is the minimizer of the strictly convex function $\sup_{\P\in\Mc(\mu)}\P[-\psi]+\nu[\psi]+ f(\psi)$, it is therefore unique. Then $\nu_l$ is unique as well. By \eqref{eq:pen_opt_duality}, $\P_{\alpha_n}$ tends to minimize $f^*(\P\circ Y^{-1}-\nu)$, by the fact that $\nu_l = \lim_{\alpha\to 0}\P_\a\circ Y^{-1}$, which concludes the proof.
\ep

\subsection{Discretization error}\label{subsect:approx_marg}

\no{\bf Proof of Proposition \ref{prop:approx_marg}} We have that $(\varphi,\psi,h)$ is a dual optimizer for $(\mu,\nu)$. Then $\varphi\oplus\psi+h^\otimes\ge c$, and if $\P'\in\Mc(\mu',\nu')$, then $\P'[c]\le \mu'[\varphi]+\nu'[\psi]\le \mu[\varphi]+\nu[\psi]+L_\varphi W_1(\mu',\mu)+L_\psi W_1(\nu',\nu)$. If we take the supremum in $\P'$, we get that
\b* 
\Sbf_{\mu',\nu'}(c)&\le& \mu[\varphi]+\nu[\psi]+L_\varphi W_1(\mu',\mu)+L_\psi W_1(\nu',\nu)\\
&=& \Ibf_{\mu,\nu}(c)+L_\varphi W_1(\mu',\mu)+L_\psi W_1(\nu',\nu)\\
&=& \Sbf_{\mu,\nu}(c)+L_\varphi W_1(\mu',\mu)+L_\psi W_1(\nu',\nu).
\e* 
As the reasoning may be symmetrical in $\big((\mu,\nu),(\mu',\nu')\big)$, we get the result.
\ep\\

\no{\bf Proof of Proposition \ref{prop:approx_marg_MC}} Similar to the proof of Proposition \ref{prop:approx_marg}, we have that
$$\Sbf_{\mu'_N,\nu'_M}(c)-\Sbf_{\mu,\nu}(c)\le (\mu_N'-\mu)[\varphi]+(\nu_M'-\nu)[\psi],$$
and
$$\Sbf_{\mu,\nu}(c)-\Sbf_{\mu'_N,\nu'_M}(c)\le (\mu-\mu_N')[\varphi_N]+(\nu-\nu_M')[\psi_M].$$
The first inequality gives
$$\Sbf_{\mu'_N,\nu'_M}(c)-\Sbf_{\mu,\nu}(c)\le (\mu_N-\mu)[\varphi]+(\nu_M-\nu)[\psi]+(\mu_N'-\mu_N)[\varphi]+(\nu_M'-\nu_M)[\psi].
$$
The two first terms are independent and their sum $(\mu_N-\mu)[\varphi]+(\nu_M-\nu)[\psi]$ is equivalent in law to $\sqrt{\frac{{\rm Var}_\mu[\varphi]}{N}+\frac{{\rm Var}_\nu[\psi]}{M}}\Nc(0,1)$ when $N,M$ go to infinity. Then doing the same work on the symmetric inequality and using the Assumptions ${\rm (i)}$ to ${\rm (iv)}$, we get the result.
\ep

\subsection{Entropy error}\label{subsect:entropy_error}

\begin{Lemma}\label{lemma:globinv}
Let $r,a>0$, $F:\R^d\longrightarrow \R^d$ and $x_0\in\R^d$ such that $\|(\nabla F)^{-1}(x_0)\|\le a^{-1}$, and on $B_{r}(x_0)$, we have that $F$ is $\Ctn^1$, that $\nabla F$ is invertible, and that $\|\nabla F-\nabla F(x_0)\|< a/2$. Then $F$ is a $\Ctn^1-$diffeomorphism on $B_{r}(x_0)$.
\end{Lemma}
\proof
We claim that $F$ is injective on $B_{r}(x_0)$, we also have that $\nabla F$ is invertible on this set. Then by the global inversion theorem, $F$ is a $\Ctn^1-$diffeomorphism on $B_{r}(x_0)$.

Now we prove the claim that $F$ is injective on $B_{r}(x_0)$. Let $x,y\in B_{r}(x_0)$,
\b*
F(y)-F(x) &=& \int_0^1\nabla F(tx+(1-t)y)(y-x)dt\\
&=& \nabla F(x)(y-x)+\int_0^1 \left[\nabla F(tx+(1-t)y)-\nabla F(x)\right](y-x)dt\\
&=& \nabla F(x) \left(y-x+\nabla F(x)^{-1}\int_0^1 \left[\nabla F(tx+(1-t)y)-\nabla F(x)\right](y-x)dt\right).
\e*
Then we assume that $F(y) = F(x)$, and we suppose for contradiction that $x\neq y$. Therefore by the fact that $\nabla F$ is invertible, we have
\b*
|y-x| &=& \left|\nabla F(x)^{-1}\int_0^1 \left[\nabla F(tx+(1-t)y)-\nabla F(x_0)+\nabla F(x_0)-\nabla F(x)\right](y-x)dt\right|\nonumber\\
&<& \|\nabla F(x)^{-1}\|2\frac{a}{2}|y-x|\nonumber\\
&\le& |y-x|,
\e*
Then we get the contradiction $|y-x|<|y-x|$. The injectivity is proved.
\ep\\

In order to prove Lemma \ref{lemma:dist_conv}, we first need the following technical lemma.

\begin{Lemma}\label{lemma:GramSchmit}
Let an integer $d\ge 1$, $k\le d$, $r>0$, and $(y_i)_{1\le i\le k+1}$, $k+1$ differentiable maps $B_r\longrightarrow\R^d$ such that $|\nabla y_i|\le A$, $|y_i|\le A$ and $\det_\aff(y_1,...,y_k)\ge A^{-1}$. Let $(u_i)_{k+1\le i\le d}$ an orthonormal basis of $\big((y_i(0)-y_{k+1}(0))_{1\le i \le k}\big)^\perp$, $(e_i)_{1\le i\le d}$ be the orthonormal basis obtained from $\big((y_i-y_{k+1})_{1\le i \le d},(u_i)_{k+1\le i \le d}\big)$ from the Gram-Schmidt process, $p$ the orthogonal projection of $0$ on $\aff(y_1,...,y_{k+1})$, and $(\lambda_i)_{1\le i\le k+1}$ the unique coefficients such that $p=\sum_{i=1}^{k+1}\lambda_i y_i$, barycentric combination. Then the maps $e_i$, $\lambda_i$, and $p$ are differentiable on $B_r$ and we may find $C,q>0$, only depending on $d$ such that if $r\le C^{-1}A^{-q}$, then we have that $|\nabla e_i|\le C A^q$, $|\nabla p|\le C A^q$, $|\nabla\lambda_i|\le C A^q$, and $\nabla \det_\aff(y_1,...,y_{k+1})\le C A^q$.
\end{Lemma}
\proof
The determinant is a polynomial expression of the coefficients, therefore if these coefficients are bounded by $A$. Then we may find $C_{\det},q_{\det}>0$ (only depending on $d$) such that $|\nabla \det|\le C_{\det}A^{q_{\det}}$.

Let $(v_i)_{1\le i \le d}:=\big((y_i-y_{k+1})_{1\le i \le d},(u_i)_{k+1\le i \le d}\big)$. Notice that for all $i$, we have $|\nabla v_i|\le 2A$, and by the fact that
\b* 
\det_\aff\big(y_1(0),...,y_{k+1}(0)\big) &:=& \left|\det\big(y_1(0)-y_{k+1}(0),...,y_k(0)-y_{k+1}(0),u_{k+1},...,u_d\big)\right|\\
&\ge& A^{-1},
\e* 
we have that $|\det(v_1,...,v_d)|\ge \frac12 A^{-1}$ on $B_{r}$ for $r\le C_{\det}^{-1}A^{-q_{\det}}\frac12 A^{-1}$.

Recall that $e_1 := \frac{v_1}{|v_1|}$. By the fact that $|v_1|...|v_d|\ge |\det(v_1,...,v_d)|$, we have that $|v_1|\ge A^{-d}$. Then we may find $C_1, q_1>0$ such that $\nabla e_1\le C_1 A^{q_1}$. Now for $1\le i\le k$, we have that $e_i := \frac{v_i-\sum_{j<i}(e_j\cdot v_i)e_j}{\left|v_i-\sum_{j<i}(e_j\cdot v_i)e_j\right|}$. Notice that
\b* 
\left|\det\left(v_1,...,v_{i-1},v_i-\sum_{j<i}(e_j\cdot v_i)e_j,v_{i+1},...,v_d\right)\right| &=& |\det\left(v_1,...,v_{i-1},v_i,v_{i+1},...,v_d\right)|\\
&\ge& A^{-1},
\e* 
and therefore we have that $\left|v_i-\sum_{j<i}(e_j\cdot v_i)e_j\right|\ge A^{-d}$. Therefore, by induction, we may find $C_i,q_i$ such that $|\nabla e_i|\le C_i A^{q_i}$.

Now notice that $p:=y_{k+1}+\sum_{i=1}^k e_i\cdot(0-y_{k+1}) e_i$. Then we may find $C_0,q_0>0$ such that $|\nabla p|\le C_0 A^{q_0}$.

Finally let $\lambda := (\lambda_1,...,\lambda_{k+1})$, $M':= \big[e_i\cdot (y_j-y_{k+1})\big]_{1\le i\le k,1\le j\le k+1}$, $M:=\left[
\begin{array}{c}
1 ... 1 \\
\hline
M'
\end{array}
\right]$, and $P:=\left[
\begin{array}{c}
1 \\
\hline
p-y_{k+1}
\end{array}
\right]$. We have that $M\lambda = P$, and therefore $\lambda = M^{-1}P$. Recall that $M^{-1} = \det(M)^{-1}Com(M)^t$ (see \eqref{eq:comatrix}), therefore we may find $C',q'>0$ such that $|M^{-1}|\le C'A^{q'}$, and $|\nabla (M^{-1})|\le C'A^{q'}$. Then we may find $C'',q''>0$ such that $|\nabla \lambda_i|\le C''A^{q''}$ for all $i$.

Finally, by the fact that
$$
\det_\aff(y_1,...,y_{k+1}) = |\det(y_1-y_{k+1},...,y_k-y_{k+1},e_{k+1},...,e_d)|,$$
we may find $C''',q'''$ such that $\det_\aff(y_1,...,y_{k+1})\le C'''A^{q'''}$.

The lemma is proved for
$$C:=\max(C_0,...,C_d,C',C'',C''',2C_{\det}),\mbox{ and for }q:=\max(q_0,...,q_d,q',q'',q''',q_{\det}+1).$$
\ep

\begin{Lemma}\label{lemma:dist_conv}
Let $A,r,e,\delta, h,H>0$ and $F:\R^d\longrightarrow \R$ such that we may find $k\in\N$ and $S\in (B_{A})^k$ such that for all $y\in S$, we have on $B_r(y)$ that $F$ is $\Ctn^2$, $A^{-1} I_d\le D^2F \le A I_d$, and $|D^2F-D^2F(y)|\le e$. Furthermore, $\det_\aff S \ge A^{-1}$, $\nabla F(S)=\{0\}$, and $\left|\sum_{i=1}^k \lambda_i S_i\right|\le H$, convex combination with $\min\lambda \ge A^{-1}$. Furthermore, assume that $F(S)\subset [-h,h]$, and $F\ge\delta\dist(Y,S)$ on $B_r(S)^c$. Then we may find $C,q>0$ such that if $\delta,r \ge CA^{q}H$, $e,H\le C^{-1}A^{-q}$, and $h\le rH$, then we have that $\big(0,F_{conv}(0)\big) = \sum_{i = 1}^{k}\overline{\lambda}_i \big(\yb_i,F(\yb_i)\big)$, convex combination with $|\yb_i-S_i|\le CA^q H$ for all $i$.
\end{Lemma}
\proof

\no\underline{Step 1:} For all $i$, the map $y\longmapsto \nabla F(y)$ is a $\Ctn^1-$diffeomorphism on $B_r(y_i)$ by Lemma \ref{lemma:globinv}. Then we define the map $z_i(a) := \nabla F^{-1}(a)$ which is defined on $B_{rA^{-1}}$. Notice that its gradient is given by $\nabla z_i(a) := D^2F^{-1}(a)$. Now we define the map $\Phi:\R^d\longrightarrow\R^d$ as follows: $\Phi_{i}(a):=a\cdot \big(z_{i}(a)- z_{k+1}(a)\big)-\Big(F\big(z_{i}(a)\big)-F\big(z_{k+1}(a)\big)\Big)$, for $1\le i\le k$, and $\Phi_i(a):=e_i(a)\cdot p(a)$ for $k+1\le i\le d$, where $p(a)$ is the orthogonal projection of $0$ on $\aff\big(z_1(a),...,z_k(a)\big)$, and $\big(e_{k+1}(a),...,e_d(a)\big)$ is the orthonormal basis of $\left(\aff\big(z_1(a),...,z_k(a)\big)\right)^\perp$ defined as the Gram-Schmidt basis obtained from $(z_1(a)-z_{k+1}(a),...,z_k(a)-z_{k+1}(a),u_{k+1},...,u_d)$, where $(u_{k+1},...,u_d)$ is a fixed basis of $\big(z_1(0)-z_{k+1}(0),...,z_k(0)-z_{k+1}(0)\big)^\perp$.

\no\underline{Step 2:} Now we prove that the convex hull $\big(F\big)_{conv}(0)$ is determined by the equation $\Phi(a)=0$ for $a$ small enough. Let $|a|\le rA^{-1}$ such that $\Phi(a) = 0$. Then $a\big(z_{i}(a)-z_{k+1}(a)\big)-\Big(F\big(z_{1}(a)\big)-F\big(z_{k+1}(a)\big)\Big) = 0$, and therefore let $b:=F\big(z_{1}(a)\big) -a z_{1}(a)=...=F\big(z_{k}(a)\big) -a z_{k}(a)$. Then the map $y\mapsto ay+b$ is tangent to $F$ at all $z_i$. Furthermore, $p(a)$ is orthogonal to $\left(\aff\big(z_1(a),...,z_k(a)\big)\right)^\perp$, implying that $p(a) = 0$. Then $0\in\aff\big(z_1(a),...,z_k(a)\big)$. By Lemma \ref{lemma:GramSchmit}, we may find $C_1,q_1>0$ (only depending on $d$) such that $|\nabla\det_\aff(z_1,...,z_k)|\le C_1A^{q_1}$. Therefore, if $a\le \frac12 A^{-2}C_1^{-1}A^{-q_1}$, we have that $|\det_\aff(z_1,...,z_k)|\ge \frac12 A^{-1}$ and we may find $(\overline{\lambda}_i)_{1\le i \le k+1}$ such that $p(a) = \sum_{i=1}^{k+1}\overline{\lambda}_i z_i(a)$, and $\overline{\lambda}_i\ge \frac12 A^{-1}$ by Lemma \ref{lemma:GramSchmit} together with the fact that $\min \lambda\ge A^{-1}$. Now we prove that $F\ge aY+b$. This holds on each $B_r\big(z_i(a)\big)$ by convexity of $F$ on these balls, together with the fact that $aY+b$ is tangent to $F$. Now out of these balls, $F\ge \delta \dist(Y,S)$ by assumption. Furthermore, $|z_i(a)-z_i(0)|\le A|a|$, and $\big|\nabla F\big(z_i(a)\big)\big|\le A^2|a|$, while similar, we have $\big|F\big(z_i(a)\big)\big|\le h+\frac12 A^3|a|^2$. Notice that as it is tangent, $aY+b = \nabla F\big(z_i(a)\big)\big(Y-z_i(a)\big)+F\big(z_i(a)\big)=\nabla F\big(z_i(a)\big)\big(Y-S_i\big)+\nabla F\big(z_i(a)\big)\big(S_i-z_i(a)\big)+F\big(z_i(a)\big)$, for all $i$. Then,
\b* 
|aY+b|&\le& \left|\nabla F\big(z_i(a)\big)\right|\left|Y-S_i\right|+\left|\nabla F\big(z_i(a)\big)\big(S_i-z_i(a)\big)+F\big(z_i(a)\big)\right|\\
&\le & A^2|a|\left|Y-S_i\right|+\frac32 A^3 |a|^2+h
\e* 
Therefore, if $\delta \ge A^2|a|+\left(\frac32 A^3 |a|^2+h\right)r^{-1}$, then $F\ge \delta\dist(Y,S)\ge aY+b$. This holds in particular if $r\ge h/H$, implying that $A^2|a|+\left(\frac32 A^3 |a|^2+h\right)r^{-1}\le A^2|a|+\frac32 A^3 |a|^2r^{-1}+hH/h \le A^2|a|+\frac32 A^3|a|+H$. Finally, the following domination is sufficient:
\be\label{eq:dom_needed_delta}
\delta \ge A^2|a|+\frac32 A^3|a|+H.
\ee 

\no\underline{Step 3:} Now we prove that $\Phi$ may be locally inverted. If $1\le i\le k$, we have $\nabla \Phi_i(a) = z_i(a)-z_{k+1}(a)$. If $k+1\le i\le d$, we have $\nabla \Phi_i(a) = \nabla p(a)e_i(a)+\nabla e_i(a)p(a)$.  We may rewrite the previous expression by introducing the locally smooth maps $\lambda_j(a)$ such that $p(a) =:\sum_{j = 1}^k\lambda_j(a)z_j(a)$, convex combination. Then $\nabla p(a) = \sum_{j = 1}^k\lambda_j(a)\nabla z_j(a)+\sum_{j = 1}^k\nabla \lambda_j(a) z_j(a)$. Notice that by the relationship $\sum_{j = 1}^k\lambda_j(a) = 1$, we have that $\sum_{j = 1}^k\nabla\lambda_j(a) = 0$, therefore $\sum_{j = 1}^k\nabla \lambda_j(a) z_j(a) = \sum_{j = 1}^k\nabla \lambda_j(a) \big(z_j(a)-z_{k+1}(a)\big)$ and $\sum_{j = 1}^k\nabla \lambda_j(a) \big(z_j(a)-z_{k+1}(a)\big)e_i = 0$ as $z_j-z_{k+1}\perp e_i$. Finally, we have $\nabla p(a)e_i(a) = \sum_{j = 1}^k\lambda_j(a)D^2F^{-1}\big(z_j(a)\big)e_i(a)$.

\no\underline{Step 4:} Now we provide a bound for $\nabla e_i(a)p(a)$. We have the control $|p(a)|\le |p(0)|+\sup_{B_{a}}|\nabla p| r\le  H +C_2A^{q_2}a$ for $a\in B_{r}$ by Lemma \ref{lemma:GramSchmit} for some $C_2,q_2> 0$. Therefore by Lemma \ref{lemma:GramSchmit}, we may find $C_3,q_3>0$ so that if $ H\le C_3^{-1} A^{-q_3}$, we have the control $|\nabla e_i(a)|\le C_3 A^{q_3}$, whence the inequality
\be\label{eq:wololo}
|\nabla e_i(a)p(a)|\le C_3 A^{q_3}( H +C_2A^{q_2}a).
\ee 

\no\underline{Step 5:} Now we provide a lower bound to $\det\nabla \Phi$. Notice that $\nabla \Phi = P_0+P'$ with $P_0:=(z_i-z_{k+1}:i\le k, De_i: i\ge k+1)$, and $P':=(\nabla e_i p: i\ge k+1)$, where $D:=\sum_{j = 1}^k\lambda_jD^2F^{-1}(z_j)$. Let $M_{basis}:= Mat(z_i-z_{k+1}:i\le k, e_i: i\ge k+1)$, then $P_0 M_{basis}^{-1}$ may be written as a block matrix as follows: $P_0 M_{basis}^{-1}=
\left[
\begin{array}{c|c}
I_k & \cdot \\
\hline
0 & D_{basis}
\end{array}
\right]$, with $D_{basis} := (e_i^tDe_j)_{k+1\le i,j\le d}$. Then $\det(P_0 M_{basis}^{-1}) = \det(D_{basis})\ge A^{-(d-k)}$, $\det(M_{basis}) = \det(z_i-z_{k+1}:i\le k)\ge A^{-1}$, and therefore $\det P_0 \ge A^{-(d-k+1)}$. Then by Lemma \ref{lemma:GramSchmit}, as $k$ lines of $P_0$ are dominated by $2A$ and $d-k$ are dominated by $A$, we have $\det\nabla\Phi(a) \ge C_4^{-1} A^{-q_4}$ if $a, H\le C_4^{-1} A^{-q_4}$, and $a\le r$ for some $C_4,q_4>0$.

\no\underline{Step 6:} Finally $|\Phi|\le A+A=2A$. In order to apply Lemma \ref{lemma:globinv}, we need to control $|\nabla\Phi(a)-\nabla\Phi(a')|$ for $a,a'\in B_r$. For $i\le k$, $|\nabla \Phi_i(a)-\nabla \Phi_i(a')| = |D^2F\big(z_i(a)\big)^{-1}-D^2F\big(z_{k+1}(a')\big)^{-1}|\le 2A^2e$. For $i\ge k+1$,
\b* 
|\nabla \Phi_i(a)-\nabla \Phi_i(a')| &=& |\nabla\big(p(a)\cdot e_i(a)\big)-\nabla\big(p(a')\cdot e_i(a')\big)|\\
&\le& |\nabla p(a)e_i(a)-\nabla p(a')e_i(a')|+|\nabla e_i(a)p(a)-\nabla e_i(a')p(a')|\\
&\le&|\nabla p(a)e_i(a)-\nabla p(a')e_i(a')|+2\big( H+C_2A^{q_2}(|a|+|a'|)\big)C_1 A^{q_1}.
\e* 
We consider the first term:
\b*
|\nabla p(a)e_i(a)-\nabla p(a')e_i(a')| &\le& \left|\sum_{j = 1}^k\lambda_j(a')\big(D^2F^{-1}(z_j(a))-D^2F^{-1}(z_j(a'))\big)\right|\\
&&+\left|\sum_{j = 1}^k\big(\lambda_j(a)-\lambda_j(a')\big)D^2F^{-1}\big(z_j(a)\big)\right|\\
&\le& A^{2}e+C_1A^{q_1}|a-a'|A,
\e*

and therefore we may find $C_5,q_5>0$ such that if $|a|,e\le C_5^{-1}A^{-q_5}$, then we have that $|\nabla \Phi(a)-\nabla \Phi(a')|\le \frac12|det \nabla \Phi(0)|^{-1}\|Com\big(\nabla \Phi(0)\big)^t\|=\|\nabla \Phi(0)\|$. Then we may apply Lemma \ref{lemma:globinv}: $\Phi$ is a $\Ctn^1-$diffeomorphism on $B_{r}$, we may find $C_6,q_6>0$ such that $C_6^{-1} A^{-q_6}\le|\nabla \Phi|\le C_6 A^{q_6}$. By assumption, we have $|\Phi(0)|\le d H$. Furthermore, $B_{rC_6^{-1} A^{-q_6}}\big(\Phi(0)\big)\subset \Phi(B_r)$. Therefore, if $ H\le rC_6^{-1}d^{-1} A^{-q_6}$, then we may find $a_0\in B_r$ such that $\Phi(a_0) = 0$. We have
\b* 
|a_0|=|\Phi^{-1}(0)-\Phi^{-1}(\Phi(0))|\le C_6A^{q_6}|\Phi(0)|\le C_6 d A^{q_6} H.
\e* 
By Step 2, $z_1(a_0),...,z_k(a_0)$ have the required property. Moreover,
$$|z_i(a_0)-S_i| = |z_i(a_0)-z_i(0)|\le C_6 d A^{q_6+1} H.$$
Finally, \eqref{eq:dom_needed_delta} is satisfied if $\delta\ge C_7A^{q_7}H$, with $C_7:=C_6+\frac52$, and $q_7 := \max(3,q_6)$. The lemma is proved for $C:= \max(3, C_1,...,C_7,C_6d)$ and $q:=\max(3,q_1,...,q_7,q_6+1)$.
\ep\\

\no{\bf Proof of Theorem \ref{thm:entropy_error}}

\no\underline{\rm Step 1:} We claim that we may find $C_1>0$ such that for $\eps$ small enough, we have $\Delta_\eps\ge -C_1\eps\ln(\eps^{-1})$, $m_\eps-$a.e. Indeed, by the fact that $(\varphi_\eps,\psi_\eps,h_\eps)$ is an optimum, we have that $e^{-\frac{\Delta_\eps}{\eps}}m_\eps$ is a probability distribution. Therefore, $e^{-\frac{\Delta_\eps(X,\cdot)}{\eps}}(m_{\eps})_{X}/\frac{d\mu_\eps}{dm_\eps^{\Xc}}$ is a probability measure, $\mu_\eps-$a.s. Then by (i), $m_\eps-$a.s., we have that
\b* 
1\ge \int_{B_{\alpha_\eps}(Y)}e^{-\frac{\Delta_\eps(X,y)}{\eps}}(m_{\eps})_{X}(dy)/\frac{d\mu_\eps}{dm_\eps^{\Xc}}&\ge& \eps^{\gamma}\int_{B_{\alpha_\eps}(Y)}e^{-\frac{\Delta_\eps(X,y)-\gamma\eps\ln(\eps^{-1})}{\eps}}(m_{\eps})_{X}(dy)\\
&\ge& \eps^{3\gamma}e^{-\frac{\Delta_\eps(X,y)}{\eps}}.
\e*

Hence $\Delta_\eps(X,Y)\ge -3\gamma\eps\ln(\eps^{-1})$, $m_\eps-$a.s. The claim is proved for $\eps$ small enough.

\no\underline{\rm Step 2:} We claim that we may find $C_2>0$ such that $\int_{\Delta_\eps\ge C_2\eps\ln(\eps^{-1})}\Delta_\eps e^{-\frac{\Delta_\eps}{\eps}}dm_\eps\ll \eps$, and $\int_{\Delta_\eps\ge C_2\eps\ln(\eps^{-1})} e^{-\frac{\Delta_\eps}{\eps}}dm_\eps\ll 1$. Indeed, let $C_2>0$. By (i), we have
$$\int_{\Delta_\eps\ge C_2\eps\ln(\eps^{-1})}\Delta_\eps e^{-\frac{\Delta_\eps}{\eps}}dm_\eps\le m_\eps[\Omega]C_2\eps\ln(\eps^{-1})\eps^{C_2}\le C_2\eps\ln(\eps^{-1})\eps^{C_2+\gamma}.$$

Similar, we have that $\int_{\Delta_\eps\ge C_2\eps\ln(\eps^{-1})} e^{-\frac{\Delta_\eps}{\eps}}dm_\eps\le C_2\eps^{C_2+\gamma}$. Therefore, up to choosing $C_2$ large enough, the claim holds.

\no\underline{\rm Step 3:} Let $\bar\Delta_\eps:=\Delta_\eps-\big(\Delta_\eps(X,\cdot)\big)_{conv}(X)-\nabla\big(\Delta_\eps(X,\cdot)\big)_{conv}(X)\cdot(Y-X)$. We claim that $\int_{x\notin D_\eps^{\Xc}}\bar\Delta_\eps e^{-\frac{\Delta_\eps}{\eps}}dm_\eps\ll \eps$. Indeed by (iii) we have that $\mu_\eps[D_\eps^{\Xc}]\ll 1/\ln(\eps^{-1})$, therefore we have
\b* 
\int_{x\notin D_\eps^{\Xc}}\Delta_\eps e^{-\frac{\Delta_\eps}{\eps}}dm_\eps
&\le&  \int_{x\notin D_\eps^{\Xc}}C_2\eps\ln(\eps^{-1}) e^{-\frac{\Delta_\eps}{\eps}}dm_\eps+\int_{\Delta_\eps\ge C_2\eps\ln(\eps^{-1})}\Delta_\eps e^{-\frac{\Delta_\eps}{\eps}}dm_\eps\\
&\le& C_2\eps\ln(\eps^{-1})\mu_\eps[(D_\eps^{\Xc})^c]+\int_{\Delta_\eps\ge C_2\eps\ln(\eps^{-1})}\Delta_\eps e^{-\frac{\Delta_\eps}{\eps}}dm_\eps\\
&\ll& \eps.
\e* 
Finally, by the martingale property of $e^{-\frac{\Delta_\eps}{\eps}}dm_\eps$, we have
\b* 
\int_{x\notin D_\eps^{\Xc}}\bar\Delta_\eps e^{-\frac{\Delta_\eps}{\eps}}dm_\eps
&=&  \int_{x\notin D_\eps^{\Xc}}\Big(\Delta_\eps-\big(\Delta_\eps(X,\cdot)\big)_{conv}(X)\\
&&-\nabla\big(\Delta_\eps(X,\cdot)\big)_{conv}(X)\cdot(Y-X)\Big) e^{-\frac{\Delta_\eps}{\eps}}dm_\eps\\
&=& \int_{x\notin D_\eps^{\Xc}}\Big(\Delta_\eps-\big(\Delta_\eps(X,\cdot)\big)_{conv}(X)\Big) e^{-\frac{\Delta_\eps}{\eps}}dm_\eps\\
&\le &\int_{x\notin D_\eps^{\Xc}}\Delta_\eps e^{-\frac{\Delta_\eps}{\eps}}dm_\eps+\mu_\eps[D_\eps^{\Xc}]C_1\eps\ln(\eps^{-1})\\
&\ll& \eps.
\e* 

\no\underline{\rm Step 4:} Let $x\notin D_\eps^\Xc$, we denote $S_x^\eps = \{s_1,...,s_{k}\}$, where $k:=k_x^\eps$. We claim that for all $S':=(s'_1,...,s'_k)\in\R^d$ such that $s'_i\in B_{r_\eps}(s_i)$ for all $i$, and $\sum \lambda'_i s'_i=x'\in B_{r_\eps}(x)$ we have for $\eps>0$ small enough that $S'\in (B_{2A_\eps})^{k}$, $|\det_\aff S'|\ge \frac12 A_\eps^{-1}$, $\min\lambda_x^\eps\ge\frac12 A_\eps^{-1}$.

Indeed, $\sum \lambda'_i (s'_i+x-x')=x$, with $s'_i+x-x'\in B_{2r_\eps}(s_i)$. By Lemma \ref{lemma:GramSchmit}, we may find $C_3,q_3$ such that $|\lambda'_i-\lambda_i|\le C_3 A_\eps^{q_3} r_\eps$, $|\det_\aff S'-\det_\aff S|\le C_3 A_\eps^{q_3} r_\eps$, and $|s'_i|\le A_\eps + r_\eps$. Now by the fact that $r_\eps \ll A_\eps^{-q_3}$, the claim is proved.

\no\underline{\rm Step 5:} We claim that up to shrinking $D_\eps^\Xc$, we may assume that for $x\notin D_\eps^\Xc$, we have $\frac{d\mu_\eps}{dm_\eps^\Xc}(x)\ge \eps^{\gamma +1}$. Indeed, $\mu_\eps^\Xc\left[\frac{d\mu_\eps}{dm_\eps^\Xc}\le \eps^{\gamma +1}\right]\le \eps^{\gamma +1}m_\eps^\Xc[\R^d]$. Therefore,
$$\mu_\eps^\Xc\left[\frac{d\mu_\eps}{dm_\eps^\Xc}\le \eps^{\gamma +1}\right]\le \eps^{\gamma +1}m_\eps^\Xc[\R^d] = \eps^{\gamma +1}m_\eps[\Omega]\le \eps\ll 1/\ln(\eps^{-1}).$$
Therefore we may shrink $D_\eps^\Xc$ by removing $\left\{\frac{d\mu_\eps}{dm_\eps^\Xc}> \eps^{\gamma+1}\right\}$ from it.

\no\underline{\rm Step 6:} We claim that for $\eps>0$ small enough, we may find unique $y_i\in B_{r_\eps}(s_i)$ such that $\nabla \Delta(x,y_i) =0$ for all $i$, $B_{\overline{r}_\eps}(y_i)\subset B_{r_\eps}(s_i)$, with $\overline{r}_\eps := \eps^{\frac12-\overline{\eta}}$, where $0<\overline{\eta}<\frac{\eta}{2}$, and finally $\Delta_\eps(x,y)\ge \frac12\sqrt{\eps}\,\dist\big(y,(y_1,...,y_k)\big)$ for $y\notin\cup_{i=1}^k B_{\bar{r}_\eps}(y_i)$. Indeed $\Delta_\eps(x,\cdot)$ is strictly convex on $B_{r_\eps}(s_i)$ for all $i$. Furthermore, let
$$m_i:=\frac{1}{\bar{\lambda}_i}\int_{B_{r_\eps}(s_i)}ye^{-\frac{\Delta_\eps(x,y)}{\eps}}dm^\Yc_\eps\left(\frac{d\mu_\eps}{dm_\eps^\Xc}\right)^{-1},$$
with $\bar{\lambda}_i :=\int_{B_{r_\eps}(s_i)}e^{-\frac{\Delta_\eps(x,y)}{\eps}}dm^\Yc_\eps\left(\frac{d\mu_\eps}{dm_\eps^\Xc}\right)^{-1}$. By the martingale property of $e^{-\frac{\Delta_\eps}{\eps}}m_\eps$, we have
$$\sum_i\bar{\lambda}_im_i = x-\int_{\left(\cup_{i}B_{r_\eps}(s_i)\right)^c}ye^{-\frac{\Delta_\eps(x,y)}{\eps}}dm^\Yc_\eps\left(\frac{d\mu_\eps}{dm_\eps^\Xc}\right)^{-1},$$ therefore by Step 5,
\b*
\left|\sum_{i=1}^k\bar{\lambda}_i m_i-x\right| &=& \left|\int_{\left(\cup_{i=1}^k B_{r_\eps}(s_i)\right)^c}y e^{-\frac{\Delta_\eps(x,y)}{\eps}}m^\Yc_\eps(dy)\left(\frac{d\mu_\eps}{dm_\eps^\Xc}\right)^{-1}\right|\\
&\le&\int_{\left(\cup_{i=1}^k B_{r_\eps}(s_i)\right)^c}|y| e^{-\eps^{-\frac12}\dist(y,S_x^\eps)}m^\Yc_\eps(dy)\eps^{-\gamma -1}\\
&\le&\int_{\left(\cup_{i=1}^k B_{r_\eps}(s_i)\right)^c, |y|\ge 2A_\eps}|y| e^{-\eps^{-\frac12}\dist(y,S_x^\eps)}m^\Yc_\eps(dy)\eps^{-\gamma -1}\\
&&+\int_{\left(\cup_{i=1}^k B_{r_\eps}(s_i)\right)^c,|y|< 2A_\eps}2A_\eps e^{-\eps^{-\frac12}\dist(y,S_x^\eps)}m^\Yc_\eps(dy)\eps^{-\gamma -1}.
\e* 
Observe that as $y_i\le A_\eps$, we have $\dist(y,S_x^\eps)\ge |y|-A_\eps$. Furthermore, if $A_\eps \ge 2$, and $|y|\ge 2A_\eps$, we have $|y|\le e^{|y|-A_\eps}$. Then we have
\be\label{eq:example_domination}
\left|\sum_{i=1}^k\bar{\lambda}_i m_i-x\right| &\le&\int_{\left(\cup_{i=1}^k B_{r_\eps}(s_i)\right)^c, |y|\ge 2A_\eps}e^{-\left(\eps^{-\frac12}-1\right)(|y|-A_\eps)}m^\Yc_\eps(dy)\eps^{-\gamma -1}\nonumber\\
&&+\int_{\left(\cup_{i=1}^k B_{r_\eps}(s_i)\right)^c,|y|< 2A_\eps}2A_\eps e^{-\eps^{-\frac12}r_\eps}m^\Yc_\eps(dy)\eps^{-\gamma -1}\nonumber\\
&\le&e^{-\left(\eps^{-\frac12}-1\right)A_\eps}m^\Yc_\eps[\R^d]\eps^{-\gamma -1}+2A_\eps e^{-\eps^{-\eta}}m^\Yc_\eps[\R^d]\eps^{-\gamma -1}\nonumber\\
&\ll&\eps.
\ee 
Similar, we have $\sum_i\bar{\lambda}_i = 1+o(\eps)$, with uniform convergence of $o(\eps)$ in $x$. By Step 4, we have that $\bar{\lambda}_i\ge \frac12 A_\eps^{-1}$ for $\eps>0$ small enough, as $\eps\ll r_\eps$. Therefore, we may find $y\in B_{r_\eps}(s_i)$ such that $\Delta_\eps(x,y)< \eps^{1-\frac{\eta}{2}}$, as otherwise, similar to \eqref{eq:example_domination}, we would have $\bar{\lambda}_i\ll\eps$. Notice that for $y$ in the boundary of $B_{r_\eps}(s_i)$, by (v), we have $\Delta_\eps(x,y)\ge \eps^{\frac12}r_\eps = \eps^{1-\eta}>\eps^{1-\frac{\eta}{2}}$, then as $\Delta_\eps(x,\cdot)$ is strictly convex on $B_{r_\eps}(s_i)$, we may find a unique minimizer $y_i\in B_{r_\eps}(s_i)$. Now let $l:=\dist\big(y_i,\partial B_{r_\eps}(s_i)\big)$, we have that $\Delta_\eps(x,y_i)+\frac12 A_\eps l^2 \ge \eps^{1-\eta}$. Then by the inequality $\Delta_\eps(x,y_i)\le \eps^{1-\frac{\eta}{2}}$, we have that $l\ge \sqrt{2A_\eps^{-1}}\eps^{\frac12-\frac{\eta}{2}}\sqrt{1-\eps^{\frac{\eta}{2}}}$. Now, let $0<\overline{\eta}<\frac{\eta}{2}$, we have $l\le \eps^{\frac12-\overline{\eta}}$ for $\eps>0$ small enough. Finally, let $y\notin\cup_{i=1}^k B_{\bar{r}_\eps}(y_i)$, we treat two cases.
\no\underline{\rm Case 1:} $y\in B_{r_\eps}(y_i)\setminus B_{\bar{r}_\eps}(y_i)$ for some $i$. Then by (iv) we have $\Delta_\eps(x,y)\ge \Delta_\eps(x,y_i)+\frac12 A_\eps^{-1}\bar{r}_\eps^2\ge -C_1\eps\ln(\eps^{-1})+\sqrt{\eps}\bar{r}_\eps\ge \frac12\sqrt{\eps}|y-y_i| \ge \frac12\sqrt{\eps}\,\dist\big(y,(y_1,...,y_k)\big)$, for $\eps>0$ small enough.

\no\underline{\rm Case 2:} $y\notin B_{r_\eps}(y_i)$ for all $i$. Then let $1\le i\le k$, we have $|y-s_i|\le r_\eps$, and recall that $|y_i-s_i|\le r_\eps$. Then $|y-y_i|\le |y-s_i|+r_\eps\le 2|y-y_i|$, and therefore $\dist\big(y,(y_1,...,y_k)\big)\le 2\dist(y,S_x^\eps)$. By (v) we have $\Delta_\eps(x,y) \ge \sqrt{\eps}\,\dist(y,S_x^\eps)\ge \frac12\sqrt{\eps}\,\dist\big(y,(y_1,...,y_k)\big)$, for $\eps>0$ small enough.

The claim is proved. Now, up to changing $\eta$ to $\overline{\eta}$, the properties (i) to (vi) are still satisfied, and the properties of (iv) and (v) also hold if we replace $S_x^\eps$ by $(y_1,...,y_k)$.

\no\underline{\rm Step 7:} Let $D_\eps^{\Xc\rightarrow\Yc}:= \left\{x\in D_\eps^\Xc:B_{r_\eps}(y)\setminus D_\eps^\Yc = \emptyset,\mbox{ for some }y\in S_x^\eps\right\}$. We claim that we have $\mu_\eps\left[D_\eps^{\Xc\rightarrow\Yc}\right]\ll 1/\ln(\eps^{-1})$. Indeed, for $x\in D_\eps^\Xc$, and for all $y\in S_x^\eps$, we have $(\P_\eps)_x[B_{r_\eps}(y)]\ge A^{-1}_\eps$ by Step 6. Therefore, if for some such $y$, we have that $B_{r_\eps}(y)\subset (D_\eps^\Yc)^c$, then $A^{-1}_\eps\le (\P_\eps)_x[B_{r_\eps}(y)]\le (\P_\eps)_x[(D_\eps^\Yc)^c]$. Then, if we integrate along $\mu_\eps$ on $D_\eps^{\Xc\rightarrow\Yc}$, together with (vi) we get that
$$A^{-1}_\eps \mu_\eps[D_\eps^{\Xc\rightarrow\Yc}]\le \P_\eps[Y\notin D_\eps^\Yc] = \nu_\eps[(D_\eps^\Yc)^c]\ll A^{-1}_\eps/\ln(\eps^{-1}).$$
The claim is proved. Now up to shrinking $D_\eps^\Xc$, we may assume that $D_\eps^\Xc\cap D_\eps^{\Xc\rightarrow\Yc}=\emptyset$.

\no\underline{\rm Step 8:} We claim that up to raising $\gamma$, if $\frac{\left|f_{|B_{2R}\setminus B_{R}}\right|}{\left\|f\right\|_\infty^{2R}}\le \frac14 \eps^\beta$ and $R\ge r_\eps/\sqrt{\eps}$, then (vi) may be applied to $y_i$ (even if $y_i\notin D_\eps^\Yc$), up to replacing $(m_{\eps})_{x}[B_{R\sqrt{\eps}}(y_i)]$ by $(m_{\eps})_{x}[B_{2R\sqrt{\eps}}(y_i')]/2^d$ for some $y_i'\in B_{r_\eps}(y_i)$. Indeed, by Step 7 we may find $y_i'\in B_{r_\eps}(y_i)\cap D_\eps^\Yc$. Let $f$ have such property. Now let $\widetilde{f}$ defined by $f= \widetilde{f}$ on $B_{R/2}$, and $\widetilde{f}(y) := \left(1-\frac{2|y|}{R}\right)f\left(\frac{R y}{2|y|}\right)$. Let $L,R\ge 1$ such that $f$ is $L-$Lipschitz, then $\widetilde{f}$ is $L-$Lipschitz. Therefore, We have $\left|\int_{B_{2R}}\widetilde{f}(y)\left[\frac{d(m_{\eps})_{x}\circ\zoom_{\sqrt{\eps}}^{y_i'}}{(m_{\eps})_{x}[B_{2R}(y_i')]}-\frac{dy}{|B_{2R}|}\right]\right|\le [2R+L]^\gamma\eps^{\beta}\int_{B_{2R}}\widetilde{f}(y)\frac{dy}{|B_{2R}|}$ from (vi). Now, as $R\ge r_\eps/\sqrt{\eps}$, we have that $B_{R\sqrt{\eps}}(y_i)\subset B_{2R\sqrt{\eps}}(y_i')$. Then
$$\left|\int_{B_{R}}f(y)\left[\frac{d(m_{\eps})_{x}\circ\zoom_{\sqrt{\eps}}^{y_i}}{(m_{\eps})_{x}[B_{2R}(y_i')]}-\frac{dy}{|B_{2R}|}\right]\right|\le [2R+L]^\gamma\eps^{\beta}\int_{B_{R}}f(y)\frac{dy}{|B_{2R}|}+|\widetilde{f}-f|_\infty.$$
As we may find $y^*\in B_{R/2}$ such that $f(y^*) = \left\|f\right\|_\infty^R$, we have $f(y)\ge \left\|f\right\|_\infty^R(1-L|y-y^*|)$, and $\int_{B_{R}}f(y)\frac{dy}{|B_{2R}|}\ge |B_1|L^{-d}$. Therefore, as $|B_{2R}| = 2^d|B_R|$, we have
\be\label{eq:neo_close_Lebesgue}
\left|\int_{B_{R}}f(y)\left[\frac{d(m_{\eps})_{x}\circ\zoom_{\sqrt{\eps}}^{y_i}}{(m_{\eps})_{x}[B_{2R}(y_i')]/2^d}-\frac{dy}{|B_{R}|}\right]\right|&\le&\left( [2R+L]^\gamma\eps^{\beta}+|B_1|^{-1}L^{d}\eps^\beta\right)\int_{B_{R}}f(y)\frac{dy}{|B_{R}|}\nonumber\\
&\le& [R+L]^{2\gamma+2}\eps^{\beta}\int_{B_{R}}f(y)\frac{dy}{|B_{R}|}.
\ee 
From now we replace $\gamma$ by $\gamma':=2\gamma+2$.

\no\underline{\rm Step 9:} As a preparation for this step, we observe that (i) to (vi) are preserved if we replace $\eta$ by any $0\le \eta'\le \eta$. Then, up to lowering $\eta$, we may assume without loss of generality that $\eta<\beta/\gamma$. Therefore
\be\label{eq:beta_not_stomped}
\beta - \eta\gamma>0.
\ee 
We claim that $\sum_{i=1}^{k_x} \int_{B_{r_\eps}(y_i)} \left(\Delta_\eps(x,y)-\Delta_\eps(x,y_i)\right)(\P_\eps)_x(dy) = \frac{d}{2}\eps + o(\eps)$, where the convergence of $o(\eps)$ is uniform in $x$. Indeed, consider
\be\label{eq:sum_at_one}
\lambda_i&:=&\int_{B_{r_\eps}(y_i)}\exp\left(-\frac{\Delta_\eps(x,y)}{\eps}\right)m^{\Yc}(dy)\frac{d\mu_\eps}{dm^{\Xc}}(x)^{-1}\\
&=&\frac{d\mu_\eps}{dm^{\Xc}}(x)^{-1}\int_{B_{\eps^{-\eta}}}\exp\left(-\frac{\Delta_\eps\left(x,\zoom_{\sqrt{\eps}}^{y_i}(y)\right)}{\eps}\right)m^{\Yc}\circ\zoom_{\sqrt{\eps}}^{y_i}(dy).\nonumber
\ee
We want to compare $\lambda_i$ to
$$\lambda_i':=\frac{d\mu_\eps}{dm^{\Xc}}(x)^{-1}\int_{B_{\eps^{-\eta}}}\exp\left(-\frac{\Delta_\eps\left(x,\zoom_{\sqrt{\eps}}^{y_i}(y)\right)}{\eps}\right)dy\frac{(m_{\eps})_{x}[B_{2r_\eps}(y_i')]/2^d}{|B_{\eps^{-\eta}}|}.$$
Notice that the map $F:y\longrightarrow \frac{\Delta_\eps\left(x,\zoom_{\sqrt{\eps}}^{y_i}(y)\right)}{\eps}$ may be differentiated in
$$\nabla F = \frac{\partial_y\Delta_\eps\left(x,\zoom_{\sqrt{\eps}}^{y_i}(y)\right)}{\sqrt{\eps}},$$
which is bounded by $A_\eps{\eps^{-\eta}}$, by the fact that $\nabla F(0)=0$ and $D^2 F\le A_\eps I_d$ by (iv). Then, we observe that the map $f: y\longrightarrow\exp\left(-\frac{\Delta_\eps\left(x,\zoom_{\sqrt{\eps}}^{y_i}(y)\right)}{\eps}\right)$ satisfies that $\frac{f}{\big\|f\big\|_\infty^{\eps^{-\eta}}}$ is $A_\eps {\eps^{-\eta}}-$Lipschitz. We may apply \eqref{eq:neo_close_Lebesgue} to $f$, with $R=r_\eps/\sqrt{\eps}$, as for $y\in B_{2R}\setminus B_R$, we have
$$\left|F(y)\right|\le e^{-\frac{\Delta_\eps(x,y_i)}{\eps}-A_\eps^{-1}|y|^2}\le |F|^{2R}_\infty e^{-A_\eps\eps^{-2\eta}}\le \frac14\eps^\beta,$$
for $\eps>0$ small enough, and get that
\be\label{eq:equiv_f}
|\lambda_i-\lambda_i'|&\le&\left|\int_{B_{{\eps^{-\eta}}}}f(y)\left[\frac{(m_{\eps})_{x}\circ\zoom_{\sqrt{\eps}}^{y_i}(dy)}{(m_{\eps})_{x}[B_{2r_\eps}(y_i')]/2^d}-\frac{dy}{|B_{{\eps^{-\eta}}}|}\right]\right|(m_{\eps})_{x}[B_{2r_\eps}(y_i')]/2^d\frac{d\mu_\eps}{dm^{\Xc}}(x)^{-1}\nonumber\\
&\le&  \left[{\eps^{-\eta}}(A_\eps+1)\right]^\gamma\eps^{\beta}\int_{B_{{\eps^{-\eta}}}} f(y)\frac{dy}{|B_{{\eps^{-\eta}}}|}(m_{\eps})_{x}[B_{2r_\eps}(y_i')]/2^d\frac{d\mu_\eps}{dm^{\Xc}}(x)^{-1}\nonumber\\
&=&(A_\eps+1)^\gamma\eps^{\beta-\eta\gamma}\lambda_i'.
\ee
Similar, we claim that the map $g$, defined by
$$g:y\longrightarrow\left(\Delta_\eps\left(x,\zoom_{\sqrt{\eps}}^{y_i}(y)\right)-\Delta_\eps\big(x,y_i\big)\right)\exp\left(-\frac{\Delta_\eps\left(x,\zoom_{\sqrt{\eps}}^{y_i}(y)\right)}{\eps}\right),$$
satisfies that $\frac{g}{\big\|g\big\|_\infty^{\eps^{-\eta}}}$ is $e^{-1}A_\eps {\eps^{-\eta}}-$Lipschitz. Now, we want to compare
$$D_i := \int_{B_{{\eps^{-\eta}}}}g(y)d(m_{\eps})_{x}\circ\zoom_{\sqrt{\eps}}^{y_i}(dy)\frac{d\mu_\eps}{dm^{\Xc}}(x)^{-1},$$
with
$$D'_i:=\int_{B_{{\eps^{-\eta}}}}g(y)dy\frac{(m_{\eps})_{x}[B_{2r_\eps}(y_i)]/2^d}{|B_{\eps^\eta}|}\frac{d\mu_\eps}{dm^{\Xc}}(x)^{-1}.$$
Hence similar, by \eqref{eq:neo_close_Lebesgue}, we have
\be\label{eq:equiv_g}
|D_i-D_i'|&=&\left|\int_{B_{{\eps^{-\eta}}}}g(y)\left[\frac{d(m_{\eps})_{x}\circ\zoom_{\sqrt{\eps}}^{y_i}(dy)}{(m_{\eps})_{x}[B_{2r_\eps}(y_i')]/2^d}-\frac{dy}{|B_{{\eps^{-\eta}}}|}\right]\right|(m_{\eps})_{x}[B_{2r_\eps}(y_i')]/2^d\frac{d\mu_\eps}{dm^{\Xc}}(x)^{-1}\nonumber\\
&\le& \left[{\eps^{-\eta}}(e^{1}A_\eps+1)\right]^\gamma\eps^{\beta}\int_{B_{{\eps^{-\eta}}}}g(x)\frac{dx}{|B_{{\eps^{-\eta}}}|}(m_{\eps})_{x}[B_{2r_\eps}(y_i')]/2^d\frac{d\mu_\eps}{dm^{\Xc}}(x)^{-1}\nonumber\\
&=&(e^{-1}A_\eps+1)^\gamma\eps^{\beta-\eta\gamma}D_i'.
\ee

Now we denote ${K_i} := \frac{(m_{\eps})_{x}[B_{2r_\eps}(y_i')]/2^d}{\left|B_{\eps^{-\eta}}\right|}\frac{d\mu_\eps}{dm^{\Xc}}(x)^{-1}\exp\left(-\frac{\Delta_\eps\big(x,y_i\big)}{\eps}\right)$, so that
\b*\label{eq:sum_o_one}
\lambda_i' =  {K_i}\int_{B_{\eps^{-\eta}}}\exp\left(-\frac{\Delta_\eps\left(x,\zoom_{\sqrt{\eps}}^{y_i}(y)\right)-\Delta_\eps\big(x,y_i\big)}{\eps}\right)dy.
\e*

We now compare $\lambda_i'$ with $\lambda_i'':=K_i\int_{\R^d}e^{-\partial_y^2\Delta_\eps(x,y_i)y^2}dy$. By the formula of the Gaussian integral, we have $\lambda_i''=K_i\sqrt{2\Pi}^{d}\sqrt{\det\partial^2_y\Delta_\eps\big(x,y_i\big)}$. Similar to \eqref{eq:example_domination}, the part of the integral out of $B_{\eps^{-\eta}}$ is uniformly negligible in front of $\eps$. We assume that $\eps>0$ is small enough so that this integral is uniformly smaller than $\eps$. By (iv), we have that
\b* 
\left(\partial^2_y\Delta_\eps\big(x,y_i\big)-\eps^\eta\right)(y-y_i)^2&\le& \Delta_\eps\left(x,\zoom_{\sqrt{\eps}}^{y_i}(y)\right)-\Delta_\eps\big(x,y_i\big)\\
&\le& \left(\partial^2_y\Delta_\eps\big(x,y_i\big)+\eps^\eta\right)(y-y_i)^2.
\e*

Therefore, we have
$$K_i\sqrt{2\Pi}^{d}\sqrt{\det\left[\partial^2_y\Delta_\eps\big(x,y_i\big)-\eps^{\eta} I_d\right]}-\eps \le\lambda_i'\le K_i\sqrt{2\Pi}^{d}\sqrt{\det\left[\partial^2_y\Delta_\eps\big(x,y_i\big)-\eps^{\eta} I_d\right]}+\eps.$$
By the fact that $\eps I_d\ll\eps^{\eta} I_d\ll A_\eps^{-1}I_d\le \partial^2_y\Delta_\eps\big(x,y_i\big)$, we may find $C_4,q_4>0$ such that 
\be\label{eq:equiv_lambda}
|\lambda_i'-\lambda_i''|\le C_4 A_\eps^{q_4}\eps^\eta \lambda_i''.
\ee 
Similar, we get that the integral of $g$ may be approximated by the integral of
$$\widetilde{g}(y) := \eps\partial^2_y\Delta_\eps\big(x,y_i\big)y^2\exp\left(-\partial^2_y\Delta_\eps\big(x,y_i\big)y^2\right).$$
Let $D_i'' := K_i\int_{\R^d}\widetilde{g}(y) dy$. Similar than the previous computation, up to raising $C_4$ and $q_4$, we have
\be\label{eq:another_domination}
|D_i'-D_i''|\le C_4 A_\eps^{q_4}\eps^\eta D_i''.
\ee 
Now we compute the value of $D_i''$. By change of variables $z = \sqrt{\partial^2_y\Delta_\eps\big(x,y_i\big)}y$, where $\sqrt{A}$ applied to a symmetrical positive definite matrix denotes the only symmetrical positive definite square root of the matrix $A$, we get that $D_i'' = \eps K_i\frac{d}{2}\sqrt{\det \partial^2_y\Delta_\eps\big(x,y_i\big)}$.

We observe that from \eqref{eq:beta_not_stomped} and (i), together with \eqref{eq:equiv_f}, \eqref{eq:equiv_g}, \eqref{eq:equiv_lambda}, and \eqref{eq:another_domination}, we have for all $i$ that $\lambda_i = \lambda_i'+o(\lambda_i')$, $\lambda_i' = \lambda_i''+o(\lambda_i'')$, $D_i = D_i'+o(D_i')$, and $D_i' = D_i''+o(D_i'')$. Finally, using \eqref{eq:prop_o} and the fact that we can sum up positive $o(\cdot)$, we get
\b* 
\sum_{i=1}^{k_x} \int_{B_{r_\eps}(y_i)} \left(\Delta_\eps(x,y)-\Delta_\eps(x,y_i)\right)(\P_\eps)_x(dy) &=& \sum_{i=1}^{k}D_i\\
&=&\sum_{i=1}^{k}D_i''+o\left(\sum_{i=1}^{k}D_i'+D_i''\right)\\
&=&\sum_{i=1}^{k}\eps K_i\frac{d}{2}\sqrt{\det \partial^2_y\Delta_\eps\big(x,y_i\big)}+o\left(\sum_{i=1}^{k}D_i\right)\\
&=&\eps \frac{d}{2}\sum_{i=1}^{k}\lambda_i''+o\left(\sum_{i=1}^{k}D_i\right)\\
&=&\eps \frac{d}{2}\sum_{i=1}^{k}\lambda_i+o\left(\sum_{i=1}^{k}D_i+\eps(\lambda_i'+\lambda_i'')\right)\\
&=&\eps \frac{d}{2}+o\left(\eps+\sum_{i=1}^{k}D_i\right)\\
&=&\eps \frac{d}{2}+o(\eps),
\e* 
where all the $o(\cdot)$ are uniform in $x$, thanks to all the controls established in this step. The claim is proved.

\no\underline{\rm Step 10:} We claim that for $\eps>0$ small enough, we have
\be\label{eq:dist_aff}
\dist\big(x,\aff(y_1,...,y_k)\big)\le C_5A_\eps^{q_5}\eps^{\frac12+\beta-\eta\gamma}+C_6A_\eps^{q_6}\eps^{\frac12+\eta}+\eps.
\ee 
Indeed, let $z_i := \lambda_i^{-1}\int_{B_{r_\eps}(y_i)}y e^{-\frac{\Delta_\eps(x,y)}{\eps}}m^\Yc_\eps(dy)\left(\frac{d\mu_\eps}{dm_\eps^\Xc}\right)^{-1}$, where recall that
$$\lambda_i = \int_{B_{r_\eps}(y_i)} e^{-\frac{\Delta_\eps(x,y)}{\eps}}m^\Yc_\eps(dy)\left(\frac{d\mu_\eps}{dm_\eps^\Xc}\right)^{-1}.$$
We have that $x = \int_{\R^d}y e^{-\frac{\Delta_\eps(x,y)}{\eps}}m^\Yc_\eps(dy)\left(\frac{d\mu_\eps}{dm_\eps^\Xc}\right)^{-1}$, by the martingale property of $e^{-\frac{\Delta_\eps}{\eps}}m_\eps$. Similar to \eqref{eq:example_domination}, we have $\left|\sum_{i=1}^k\lambda_i z_i-x\right| \ll\eps$. Similar, we also have $\sum_{i=1}^k\lambda_i = 1 + o(\eps)$. Therefore, $\dist\big(x,\aff(z_1,...,z_k)\big)\ll \eps$. Now let
$$z_i':=\frac{1}{\lambda_i'}\int_{B_{r_\eps}(y_i)}y e^{-\frac{\Delta_\eps(x,y)}{\eps}}dy\frac{(m_{\eps})_{x}[B_{2r_\eps}(y_i')]/2^d}{|B_{r_\eps}|}\left(\frac{d\mu_\eps}{dm_\eps^\Xc}\right)^{-1},$$
where recall that $\lambda_i'=\int_{B_{r_\eps}(y_i)} e^{-\frac{\Delta_\eps(x,y)}{\eps}}dy\frac{(m_{\eps})_{x}[B_{r_\eps}(y_i)]}{|B_{r_\eps}|}\left(\frac{d\mu_\eps}{dm_\eps^\Xc}\right)^{-1}$.
We claim that we may find universal $C_5,q_5>0$, such that $|\lambda_i(z_i-y_i)-\lambda_i'(z_i'-y_i)|\le C_5A_\eps^{q_5}\eps^{\frac12+\beta-\eta\gamma}$. Indeed, if $u$ is a unit vector, we have
$$h:y\longrightarrow \sqrt{\eps}|y\cdot u|\exp\left(-\frac{\Delta_\eps\left(x,\zoom_{\sqrt{\eps}}^{y_i}(y)\right)-\Delta_\eps(x,y_i)}{\eps}\right).$$
We claim that $\frac{h}{\|h\|_\infty^{\eps^{-\eta}}}$ is $\eps^{-\eta}A_\eps\left(1+\frac{\sqrt{A_\eps}}{2y_0^3}\right)-$Lipschitz, where $y_0>0$ is the unique positive real satisfying $2y_0^2 = e^{-y_0^2}$, furthermore the function is small enough on $B_{2R}(y_i)\setminus B_{R}(y_i)$ so that we may apply \eqref{eq:neo_close_Lebesgue}:
\be\label{eq:Gaussian_estimate}
\left|u\cdot\big(\lambda_i(z_i-y_i)-\lambda_i'(z_i'-y_i)\big)\right|&=&\left|\int_{B_{\eps^{-\eta}}(y_i)}h(y)\left[\frac{d(m_{\eps})_{x}\circ\zoom_{\sqrt{\eps}}^{y_i}(dy)}{(m_{\eps})_{x}[B_{2r_\eps}(y_i')]/2^d}-\frac{dy}{|B_{{\eps^{-\eta}}}|}\right]\right|K_i|B_{\eps^{-\eta}}|\nonumber\\
&\le& \left[\eps^{-\eta}A_\eps\left(1+\frac{\sqrt{A_\eps}}{2y_0^3}\right)\right]^\gamma\eps^\beta\int_{B_{\eps^{-\eta}}(y_i)}h(y)dy K_i,\nonumber\\
&\le&\left[A_\eps\left(1+\frac{\sqrt{A_\eps}}{2y_0^3}\right)\right]^\gamma\eps^{\frac12+\beta-\gamma\eta} I A_\eps  K_i,
\ee 
where $I:=\int_\R^d |y\cdot u|e^{-|y|^2}dy$, as recall that by Step 9, $K_i=\frac{\lambda_i''}{\sqrt{2\Pi}^d\sqrt{\det\partial_y^2\Delta_\eps(x,y_i)}}\le A_\eps^{d/2}$. Now we consider again a unit vector $u$, and
\b* 
(z_i'-y_i)\cdot u &=& \frac{\sqrt{\eps}}{\lambda_i'}\int_{B_{\eps^{-\eta}}}\big((y\cdot u)_+-(y\cdot u)_-\big) e^{-\frac{\Delta_\eps(x,y)-\Delta_\eps(x,y_i)}{\eps}}dy K_i\\
&\le& \frac{\sqrt{\eps}}{\lambda_i'}\int_{B_{\eps^{-\eta}}}(y\cdot u)_+ e^{-(\partial_y^2\Delta_\eps(x,y_i)-\eps^\eta I_d)y^2}dy K_i\\
&&-\frac{\sqrt{\eps}}{\lambda_i'}\int_{B_{\eps^{-\eta}}}(y\cdot u)_- e^{-(\partial_y^2\Delta_\eps(x,y_i)+\eps^\eta I_d)y^2}dy K_i\\
&\le& \frac{K_i\sqrt{\eps}}{\lambda_i'}(I/2)\left(\sqrt{\partial_y^2\Delta_\eps(x,y_i)-\eps^\eta I_d}^{-1}u^2-\sqrt{\partial_y^2\Delta_\eps(x,y_i)+\eps^\eta I_d}^{-1}u^2\right)\\
&\le& C_6 A_\eps^{q_6}\eps^{\frac12+\eta},
\e*
for some $C_6,q_6>0$, independent of $x$ and $u$, as $\partial_y^2\Delta_\eps(x,y_i)\ge A_\eps^{-1}I_d\gg \eps^\eta I_d$, and $\lambda_i'\ge \frac12A_\eps^{-1}$. Then $|z_i'-y_i|\le C_6 A_\eps^{q_6}\eps^{\frac12+\eta}$.

Finally, with \eqref{eq:equiv_f} and
\eqref{eq:equiv_lambda}, up to raising $C_5,C_6,q_5,q_6$, we get the estimate $|z_i-y_i|\le C_5 A_\eps^{q_5}\eps^{\frac12+\beta-\gamma\eta}+ C_6 A_\eps^{q_6}\eps^{\frac12+\eta}$. We finally get the desired estimate from the fact that $\dist\big(x,\aff(y_1,...,y_k)\big)$.

\no\underline{\rm Step 11:} We now claim that
\be\label{eq:error_conv}
\int_{B_{r_\eps}(y_i)}\bar\Delta_\eps(x,\cdot) d(\P_\eps)_x = \int_{B_{r_\eps}(y_i)}\big(\Delta_\eps(x,\cdot) -\Delta_\eps(x,y_i) \big)d(\P_\eps)_x +o(\eps),
\ee 
where the convergence speed of $o(\eps)$ is independent of the choice of $x$ and $i$. Indeed, by \eqref{eq:dist_aff} and the fact that $$|\Delta_\eps(x,y_i)|\le \max(C_1,C_2)\eps\ln(\eps^{-1})\le \eps^{1-\frac{\eta}{2}} \le r_\eps H_\eps,$$
with $H_\eps := \eps^{\frac12+\frac12\min(\eta, \beta-\eta\gamma)}$ for $\eps$ small enough. Furthermore, by \eqref{eq:dist_aff}, we have that for $\eps>0$ small enough, $\dist\big(x,\aff(y_1,...,y_k)\big)\le H_\eps$. Then, we may apply Lemma \ref{lemma:dist_conv}: we may find $C_7,q_7>0$ such that we have $\bar\Delta_\eps(x,y) = \Delta_\eps(x,y)-\Delta_\eps(x,\yb_i)-\nabla\Delta_\eps(x,\yb_i)\cdot(y-\yb_i)$, with $|y_i-\yb_i|\le C_7 A_\eps^{q_7}H_\eps\le C_7 A_\eps^{q_7}\eps^{\frac12}$, as by Step 6, we have that for $y\notin \cup_{i=1}^k B_{r_\eps}(y_i)$, we have $\Delta_\eps(x,y)\ge \frac12\sqrt{\eps}\,\dist\big(y,(y_1,...,y_k)\big)$, where $\frac12\sqrt{\eps}\gg C_7 A_\eps^{q_7} H_\eps$, and $r_\eps\gg C_7 A_\eps^{q_7} H_\eps$ as well. Notice that therefore, we have $0\le \Delta_\eps(x,\yb_i)-\Delta_\eps(x,y_i)\le \frac12 A_\eps C_7^2 A_\eps^{2q_7}H_\eps^2\ll \eps$, and similar, we have $|\nabla \Delta_\eps(x,\yb_i)\cdot (y_i-\yb_i)|\le A_\eps C_7^2 A_\eps^{2q_7}H_\eps^2\ll \eps$. Finally, $\int_{B_{r_\eps}(y_i)}\nabla\Delta_\eps(x,\yb_i)\cdot (y-y_i) d(\P_\eps)_x = \nabla\Delta_\eps(x,\yb_i)\cdot\int_{B_{r_\eps}(y_i)}(y-y_i) d(\P_\eps)_x$. We have from the computations in Step 10 that $\left|\int_{B_{r_\eps}(y_i)}(y-y_i) d(\P_\eps)_x\right|\le C_3 A_\eps^{q_3}\eps^{\frac12}$, and $|\nabla\Delta_\eps(x,\yb_i)\cdot\int_{B_{r_\eps}(y_i)}(y-y_i) d(\P_\eps)_x|\le C_3 C_7 A_\eps^{q_3+2q_7}\eps^{\frac12}H_\eps \ll \eps$, whence the result.

\no\underline{\rm Step 12:} Now using Step 3 and Step 9, integrating against $\mu_\eps$, with the uniform error estimate \eqref{eq:error_conv}, together with controls that are independent of $x$, similar to \eqref{eq:example_domination}, to deal with $\left(\cup_{i=1}^k B_{r_\eps}(y_i)\right)^c$, we get
\b*
\int\bar\Delta_\eps d\P_\eps &=& \int_{D_\eps^\Xc} \sum_{i=1}^{k_x} \int_{B_{r_\eps}(y_i)} \left(\Delta_\eps(x,y)-\Delta_\eps(x,y_i)\right)(\P_\eps)_x(dy)\mu_\eps(dx)+o(\eps)\\
&=& \eps \frac{d}{2} +o(\eps).
\e*
Finally, notice that
\b* 
\bar\Delta_\eps&=&\Delta_\eps-\big(\Delta_\eps(X,\cdot)\big)_{conv}(X)-\nabla\big(\Delta_\eps(X,\cdot)\big)_{conv}(X)\cdot(Y-X)\\
&=& \psi_\eps - c- \big(\psi_\eps-c(X,\cdot)\big)_{conv}(X)-\nabla\big(\psi_\eps-c(X,\cdot)\big)_{conv}(X)\cdot(Y-X)\\
&=& \bar\varphi_\eps + \psi_\eps - c-\nabla\big(\psi_\eps-c(X,\cdot)\big)_{conv}(X)\cdot(Y-X).
\e* 
Therefore we have $\P_\eps[\bar\Delta_\eps] = \mu_\eps[\bar\varphi_\eps]+\nu_\eps[\psi_\eps]-\P_\eps[c]$. Whence the result of the theorem.
\ep

\subsection{Asymptotic penalization error}\label{subsect:limit_pen}

\no{\bf Proof of Proposition \ref{prop:limit_pen}} As $\mu\preceq_c\nu$, $\widetilde{V}_\eps$ is convex with a finite global minimum. Then the minimum of $\widetilde{V}_\eps+\alpha f$ converges to a minimum of $\widetilde{V}_\eps$. More precisely, let $\psi_\alpha$ be the only global minimizer of $\widetilde{V}_\eps+\alpha f$, then $\psi_\a$ is also the minimizer of the map $\frac{1}{\alpha}\left(\widetilde{V}_\eps-\widetilde{V}_\eps(\psi_0)\right)+ f$, which $\Gamma-$converges to $\iota_{\left\{\psi:\widetilde{V}_\eps(\psi) \neq \widetilde{V}_\eps(\psi_0)\right\}}+f$, whose unique global minimizer is $\psi_0$. Therefore $\psi_\a\longrightarrow\psi_0$ when $\a\longrightarrow 0$. Now the first order condition gives that $\frac{\nu_\a - \nu}{\a} = \nabla f(\psi_\a)\longrightarrow \nabla f(\psi_0)$, when $\a\longrightarrow 0$, by convexity and differentiability of $f$, guaranteeing that $\psi\longmapsto \nabla f(\psi)$ is continuous.
\ep

\subsection{Convergence of the martingale Sinkhorn algorithm}\label{subsect:cvg_Sinkhorn}




\no{\bf Proof of Theorem \ref{thm:cvg_Sinkhorn}}
This result stems from an indirect application of Theorem 5.2 in \cite{beck2013convergence}. By a direct application of this theorem we get that
\be 
&V_\eps(x_k)-V_\eps^*\le \frac{2\min(L_1,L_2)R^2(x_0)}{n-1},\label{eq:convergence_Beck1}\\
\mbox{and that}&V_\eps(x_k)-V_\eps^*\le \left(1-\frac{\sigma}{\min(L_1,L_2)}\right)^{n-1}\big(V_\eps(x_0)-V_\eps^*\big),\label{eq:convergence_Beck2}
\ee 
with $R(x_0):=\sup_{V_\eps(x)\le V_\eps(x_0)}\dist(x,\Xc^*)$, $L_1$ (resp. $L_2$) is the Lipschitz constant of the $\psi-$gradient (resp. $(\varphi,h)-$gradient) of $V_\eps$, and $\sigma$ is the strong convexity parameter of $V_\eps$. Furthermore, the strong convexity gives that
\be\label{eq:distance_from_conv}
\dist(x_k,\Xc^*)\le \sqrt{\frac{2}{\sigma}}(V_\eps(x_k)-V_\eps^*)^{\frac12}.
\ee 
However the gradient $\nabla V_\eps$ is locally but not globally Lipschitz, nor $V_\eps$ strongly convex. Therefore we need to refine the theorem by looking carefully at where these constants are used in its proof.

\no\underline{\rm Step 1:} The constant $L_1$ is used for Lemma 5.1 in \cite{beck2013convergence}. We need for all $k\ge 0$ to have $V_\eps(x_k)-V_\eps(x_{k+1/2})\ge \frac{1}{2 L_1}|\nabla V_\eps(x_k)|^2$. We may start from $x_1$, this way all the $x_k$ are such that $(e^{-\frac{\Delta(x_k)_{i,j}}{\eps}})_{i\in\Xc,j\in \Yc}$ is a probability distribution. Then $|\partial^2_\psi V_\eps(x_k)|\le \eps^{-1}$. Let $u\in \R^{\Yc}$, then $|\partial^2_\psi V_\eps(\varphi_k,\psi_k+u,h_k)|\le \eps^{-1}e^{\frac{|u|_\infty}{\eps}}$. We want to find $C,L>0$ such that $V_\eps(x_k)-V_\eps(\varphi_k,\psi_k-C\partial_\psi V_\eps(x_k),h_k)\ge \frac{1}{2L}|\partial_\psi V_\eps(x_k)|$, then $L$ may be use to replace $L_1$ in the final step of the proof of Lemma 5.1 in \cite{beck2013convergence}. Recall that $|\partial_\psi V_\eps(x_k)|_\infty\le 1$, as it is the difference of two probability vectors. We have
\b* 
V_\eps(x_k)-V_\eps\big(\varphi_k,\psi_k-C\partial_\psi V_\eps(x_k),h_k\big)\\
= -\partial_\psi V_\eps(x_k)\cdot \big(-C\partial_\psi V_\eps(x_k)\big)\\
-C^2\int_0^1 (1-t)\partial_\psi V_\eps(x_k)^t \partial^2_\psi V_\eps\big(\varphi_k,\psi_k-tC\partial_\psi V_\eps(x_k),h_k\big)\partial_\psi V_\eps(x_k) dt\\
\le C|\partial_\psi V_\eps(x_k)|^2-C^2|\partial_\psi V_\eps(x_k)|^2\int_0^1\eps^{-1}(1-t)e^{\frac{tC}{\eps}} dt\\
= C|\partial_\psi V_\eps(x_k)|^2-C^2|\partial_\psi V_\eps(x_k)|^2\eps^{-1}\frac{e^{\frac{C}{\eps}}-1-\frac{C}{\eps}}{\frac{C^2}{\eps^2}}\\
=\left(C-\eps \left(e^{\frac{C}{\eps}}-1-\frac{C}{\eps}\right)\right)|\partial_\psi V_\eps(x_k)|^2.
\e* 
Deriving with respect to $C$ gives the equation $C = \eps\ln(2)$. We get
\b* 
V_\eps(x_k)-V_\eps\big(\varphi_k,\psi_k-C\partial_\psi V_\eps(x_k),h_k\big)\ge \eps \big(2\ln(2)-1\big)|\partial_\psi V_\eps(x_k)|^2.
\e* 
Therefore we may use $L := \eps^{-1}\big(4\ln(2)-2\big)^{-1}$.

\no\underline{\rm Step 2:} The constant $\sigma$ is used to get the result from (3.21) in \cite{beck2013convergence}. Then we just need the inequality
\be\label{eq:ineq_sec_order}
V_\eps(y)\ge V_\eps(x)+\nabla V_\eps(x)\cdot(y-x)+\frac{\sigma}{2}|y-x|^2,
\ee 

to hold for some $y \in \Xc^*$ and $x = x_k$ for all $k\ge 0$. Now we give a lower bound for $\sigma$. Notice that $V_\eps = \mu[\varphi]+\nu[\psi]+ \eps\sum_{x\in\Xc,y\in\Yc}exp\left(-\frac{\cdot}{\eps}\right)\circ \Delta_{x,y}$. Then for $x_0,u\in \Dc_{\Xc,\Yc}$, we have
\b* 
u^tD^2V_\eps(x_0) u &=& \eps^{-1}\sum_{x\in\Xc,y\in\Yc}exp\left(-\frac{\cdot}{\eps}\right)\circ \Delta_{x,y}(x_0)\Delta_{x,y}(u)^2\\
&\ge& \eps^{-1}exp\left(-\frac{|\Delta(x)|_\infty}{\eps}\right)|\Delta(u)|^2.
\e*

Then, by definition of $\lambda_2$, we may find $\widetilde{u}$ such that $\Delta(u) = \Delta(\widetilde{u})$, and
\be\label{eq:coerc}
u^tD^2V_\eps(x_0) u \ge \frac{|\Xc|}{\lambda_2\eps}exp\left(-\frac{|\Delta(x)|_\infty}{\eps}\right)|\widetilde{u}|^2.
\ee 
Now, we claim that $|\Delta(x)|_\infty\le D(x_0)$. Then let $x^*\in\Xc^*$ and consider \eqref{eq:coerc} for $u = x^*-x$. Then we have that $x+\widetilde{u}\in\Xc^*$, and therefore, we may take $y = x+\widetilde{u}$ for \eqref{eq:ineq_sec_order}, and therefore use
\be\label{eq:value_sigma}
\sigma := \frac{|\Xc|}{\lambda_2\eps}exp\left(-\frac{D(x_0)}{\eps}\right),
\ee 
in this equation.

\no\underline{\rm Step 3:} Now we prove our claim that $|\Delta(x)|_\infty\le D(x_0)$. Indeed
$$V_\eps(x_0)\ge V_\eps(x)=\mu[\varphi]+\nu[\psi]+\eps=\P_0[\varphi\oplus\psi + h^\otimes]+\eps = \P_0[\Delta(x)]+\P_0[c]+\eps.$$

Therefore we have $\P_0[\Delta(x)]\le V_\eps(x_0)-\P_0[c]-\eps$, and finally
\be\label{eq:dom_norm1_Delta}
(\P_0)_{min}|\Delta(x)|_1\le V_\eps(x_0)-\P_0[c]-\eps.
\ee 
Then $|\Delta(x)|_\infty\le D(x_0)$ stems from the definition of $\lambda_1$.

\no\underline{\rm Step 4:} Now we provide the bound on $R(x_0)$. From the proof of Theorem 5.2 in \cite{beck2013convergence}, what is needed to make the proof work is $R(x_0) = \sup_{k\ge 0}\dist(x_k,\Xc^*)$, which is smaller than $\sup_{V_\eps(x)\le V_\eps(x_0)}\dist(x,\Xc^*)$. Furthermore, from \eqref{eq:distance_from_conv} together with \eqref{eq:value_sigma}, we get that the supremum $\sup_{k\ge 0}\dist(x_k,\Xc^*)$ is also smaller than $\sqrt{\frac{\lambda_{2}\eps}{|\Xc|} e^{\frac{D(x_0)}{\eps}}} \big(V_\eps(x_0)-V_\eps^*\big)^{\frac12}$. Finally, from \eqref{eq:dom_norm1_Delta} together with the definition of $\lambda_1$, we may find that $\widetilde{x}^*,\widetilde{x}_k\in \Dc_{\Xc,\Yc}$ such that $\Delta(x_k) = \Delta(\widetilde{x}_k)$, $\widetilde{x}^* \in \Xc^*$, $
|\widetilde{x}_k|_1\le \lambda_1\frac{V_\eps(x_0)-\P_0[c]-\eps}{|\Xc|(\P_0)_{\min}}
$, and $
|\widetilde{x}^*|_1\le \lambda_1\frac{V_\eps(x_0)-\P_0[c]-\eps}{|\Xc|(\P_0)_{\min}}
$. Then $|\widetilde{x}_k-\widetilde{x}^*|\le |\widetilde{x}_k-\widetilde{x}^*|_1\le 2\lambda_1\frac{V_\eps(x_0)-\P_0[c]-\eps}{|\Xc|(\P_0)_{\min}}$ by the fact that $|\cdot|\le |\cdot|_1$. Finally, as $\widetilde{x}^*+x_k-\widetilde{x}_k\in\Xc^*$, we have $\dist(x_k,\Xc^*)\le 2\lambda_1\frac{V_\eps(x_0)-\P_0[c]-\eps}{|\Xc|(\P_0)_{\min}}$. Therefore the same bound holds for $\sup_{k\ge 0}\dist(x_k,\Xc^*)$.

\no\underline{\rm Step 5:} Finally, as we focus on the $L_1$ optimization phase, we may replace $n-1$ by $n$ in the convergence formula \eqref{eq:convergence_Beck1} and \eqref{eq:convergence_Beck2}, see the proof of Theorem 5.2 in \cite{beck2013convergence}. The result is proved.
\ep

\subsection{Implied Newton equivalence}\label{subsect:equivNewton}

\no{\bf Proof of Proposition \ref{prop:equivNewton}} We apply the Newton step in the algorithm to $\big(x,y(x)\big)$. we are looking for $p$ such that $D^2Fp = \nabla F$. First $\nabla F\big(x,y(x)\big) = \partial_xF\big(x,y(x)\big) = \nabla \tilde{F}(x)$, then if we decompose $p = p_x\oplus p_y\in\Xc\oplus\Yc$, the equation becomes
\b*
\partial^2_xFp_x+\partial^2_{xy}Fp_y = \nabla \tilde{F}(x),&\mbox{and}&\partial^2_{yx}Fp_x+\partial^2_yFp_y =0.
\e*
The solution to this equation system is given by $p_y = -\partial^2_yF^{-1}\partial^2_{yx}Fp_x$, and
\be\label{eq:Newtonimplied}
\big(\partial^2_xFp_x-\partial^2_{xy}F\partial^2_yF^{-1}\partial^2_{yx}F\big)p_x = \nabla \tilde{F}(x)
\ee
The conclusion follows from the fact that \eqref{eq:Newtonimplied} is the step for the Newton algorithm applied to $\tilde{F}$. The Newton step on $y$ does not matter, as $y$ will be immediately thrown away and replaced by $y(x)$.
\ep

\section{Numerical experiment}\label{sect:numerics}

\subsection{An hybrid algorithm}

The steps of the Newton algorithm are theoretically very performing if the current point is close enough to the optimum. What is really time-consuming is the computation of the descent direction with the conjugate gradient algorithm. The idea of preferring the Newton method to the Bregman projection method in the case of martingale optimal transport comes from the fact that, unlike in the case of classical transport, projecting on the martingale constraint is more costly than projecting on the marginal constraints, as we use a Newton algorithm instead of a closed formula. From the experiment, we would say that in dimension $1$ it takes $7$ times more time, and $20$ times more in dimension $2$. The implied Newton algorithm performs this projection only for the Newton step, whereas it is not necessary for the conjugate gradient algorithm.

We Notice that the Bregman projection algorithm is more effective at the beginning, to find the optimal region, and then it converges slower. In contrast, the Newton algorithm is slow at the beginning when it is searching the neighborhood of the optimum, but when its finds this neighborhood, the convergence gets very fast. Then it makes sense to apply an hybrid algorithm that starts with Bregman projections, and concludes with the Newton method. We call this dual-method algorithm the hybrid algorithm. We see on the simulations that it generally out-performs the two other algorithms.

\begin{figure}[h]
\centering
\begin{subfigure}{.49\textwidth}
  \centering
  \includegraphics[width=.99\linewidth]{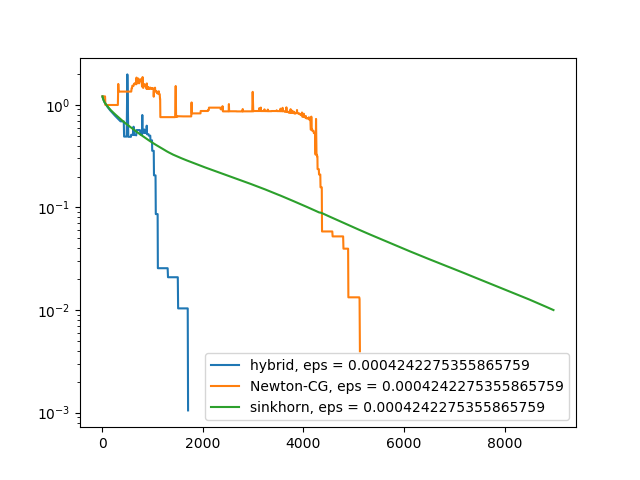}
  \caption{Dimension 1.}
  \label{fig:performancesdim1}
\end{subfigure}%
\begin{subfigure}{.49\textwidth}
  \centering
  \includegraphics[width=.99\linewidth]{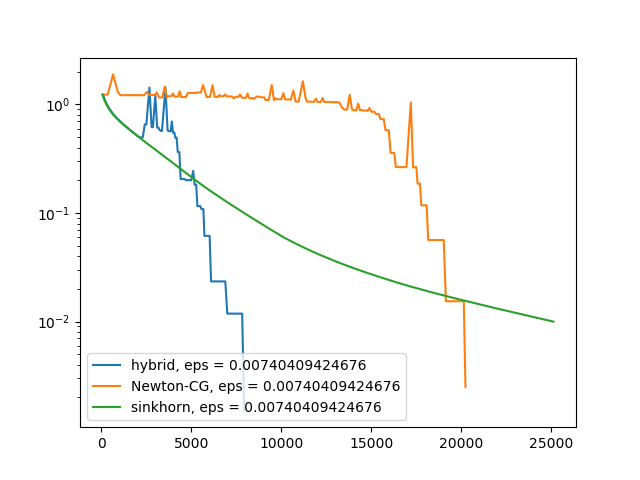}
  \caption{Dimension 2.}
  \label{fig:performancesdim2}
\end{subfigure}
\caption{Log plot of the size of the gradient VS time for the Bregman projection algorithm, the Newton algorithm, and the Hybrid algorithm.}
\label{fig:performances}
\end{figure}

Figure \ref{fig:performances} compares the evolution of the gradient error in dimension $1$ and $2$ of the longest step of the three algorithms in terms of computation time. What we call here the gradient error is the norm $1$ of the gradient of the function $\widetilde{V}_\eps$ that we are minimizing, and which is also equal to the difference between the target measure $\nu$ and the current measure. In the case of Newton algorithms, the penalization gradient is also included, then we use a coefficient in front of this penalization so that it does not interfere too much with the equation between the current and the target measure. We use the $\eps-$scaling technique. For each value of $\eps$, we iterate the minimization algorithm until the error is smaller than $10^{-2}$. Then at the final iteration we lower the target error to the one we want.

The green line corresponds to the Bregman projections algorithm. The orange line corresponds to the implied truncated Newton algorithm. All the techniques evocated in Section \ref{sect:practi} are applied. We use the diagonal of the Hessian to precondition the conjugate gradient algorithm. The coefficient in front of the quadratic penalization, which is normalized by $\nu^2$, is set to $10^{-2}$. Finally the blue line corresponds to the "hybrid algorithm", which consists in doing some Bregman projection steps before switching to the implied truncated Newton algorithm. The moment of switching is chosen by very empirical criteria: we do it after having the initial error divided by $2$ or after $100$ iteration, or if the initial error is divided by $1.1$ if the initial error is smaller than $0.1$.

Figure \ref{fig:performancesdim1} gives the computation times of these three entropic algorithms, for a grid size going from $10$ to $2500$ while $\eps$ goes from $1$ to $10^{-4}$, with the cost function $c := XY^2$, $\mu$ uniform on $[-1,1]$, and $\nu := \frac{1}{K}|Y|^{1.5}\mu$, where $K := (|Y|^{1.5}\mu)[\R]$. By \cite{henry2016explicit} the optimal coupling that we get is the "left curtain" coupling studied in \cite{beiglboeck2016problem}. We show the curves for the value of $\eps$ that takes the largest amount of time, the one for which the time of computation is the most important for $\eps = 4.2\x 10^{-4}$.

We conduct the same experiment on a two dimensional problem. The difference of efficiency between the algorithms should be even bigger, as the computing of the optimal $h$ becomes more costly, as the optimization of a convex function of two variables. Let $d=2$, $c:(x,y)\in\R^2\x\R^2\longmapsto x_1(y_1^2+2y_2^2)+x_2(2y_1^2+y_2^2)$, $\mu$ uniform on $[-1,1]^2$, and $\nu = \frac{1}{K}(|Y_1|^{1.5}+|Y_2|^{1.5})\mu$ where $K:= \big((|Y_1|^{1.5}+|Y_2|^{1.5})\mu\big)[\R^2]$. We start with a $10\x10$ grid and scale it to a $160\x160$ one while $\eps$ scales from $1$ to $10^{-4}$. Figure \ref{fig:performancesdim2} gives the computation times of the three entropic algorithms. Once again we show the curves for the value of $\eps$ that takes the largest amount of time, the one for which the time of computation is the most important for $\eps = 7.4\x10^{-3}$.

\subsection{Results for some simple cost functions}

\subsubsection{Examples in one dimension}

\begin{figure}[h]
\centering
\begin{subfigure}{.32\textwidth}
  \centering
  \includegraphics[width=.99\linewidth]{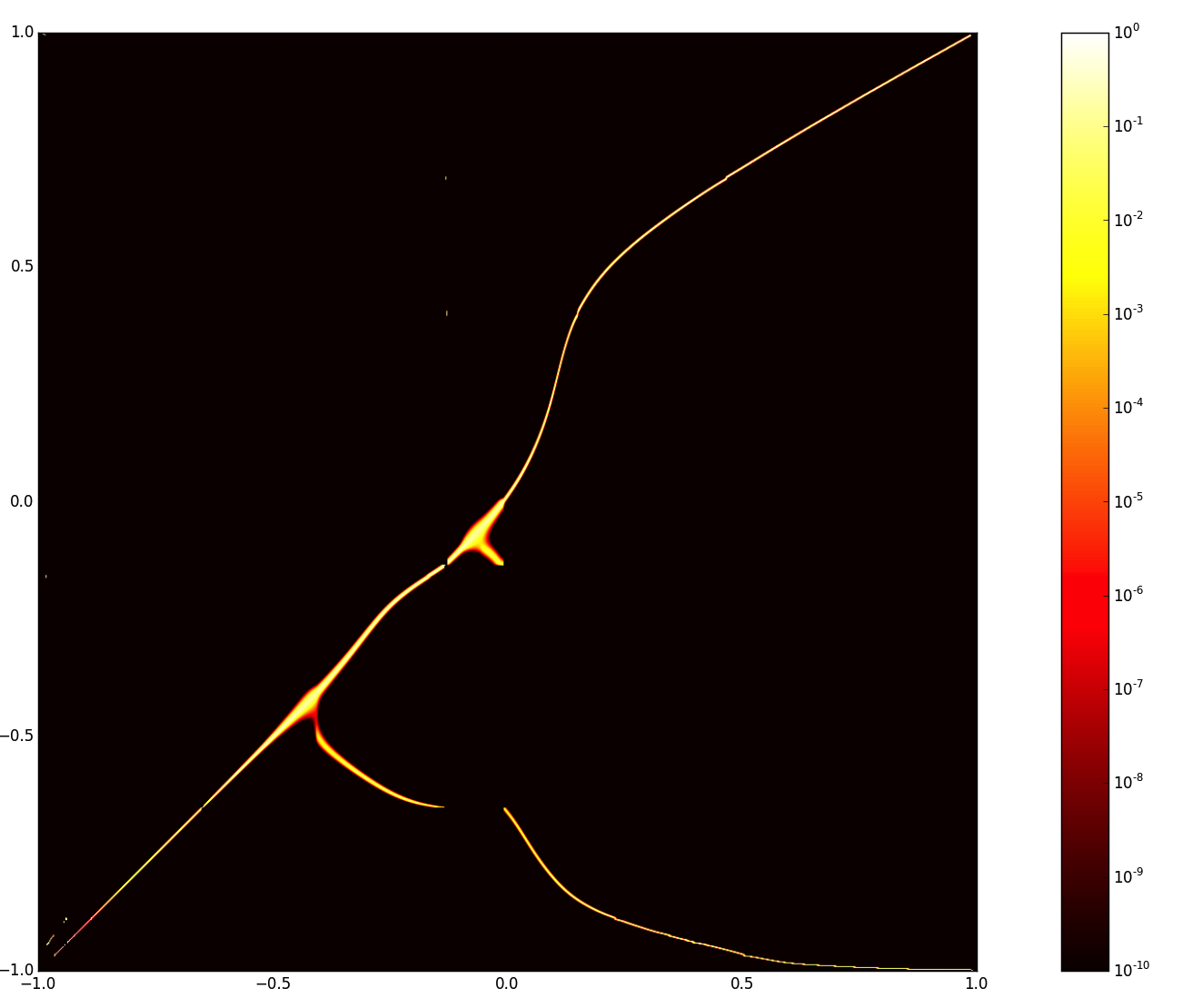}
  \caption{$c = XY^2$.}
  \label{fig:left_curtain}
\end{subfigure}%
\begin{subfigure}{.32\textwidth}
  \centering
  \includegraphics[width=.99\linewidth]{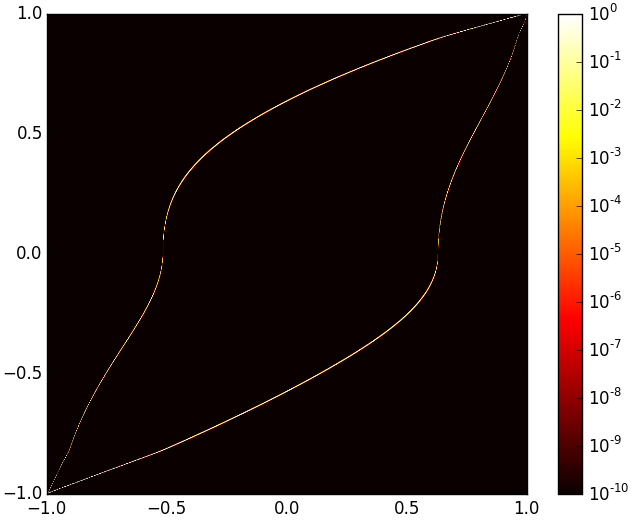}
  \caption{$c = |X-Y|$.}
  \label{fig:distance_cost}
\end{subfigure}
\begin{subfigure}{.32\textwidth}
  \centering
  \includegraphics[width=.99\linewidth]{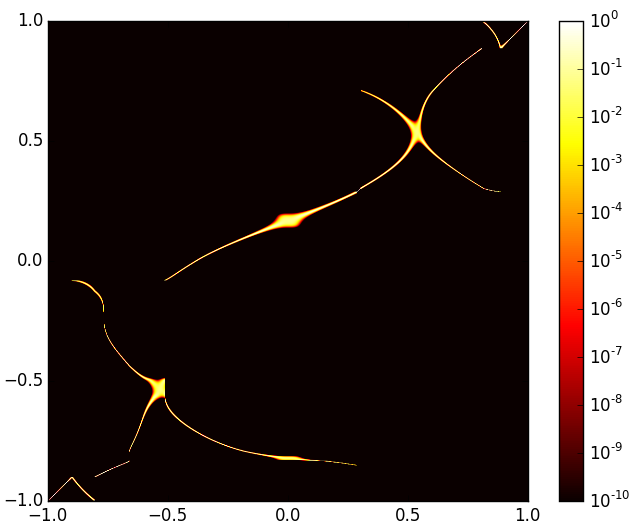}
  \caption{$c = \sin(8XY)$.}
  \label{fig:random_cost}
\end{subfigure}
\caption{Optimal coupling for different costs in dimension one.}
\label{fig:solution_dim1}
\end{figure}

Figure \ref{fig:solution_dim1} give the solution for three different costs for $\eps = 10^{-5}$ with $\mu := (\mu_1+\mu_2)/2$ and $\nu := (\nu_1+\nu_2)/2$ with $mu_1$ uniform on $[-1,1]$, $\nu_1 = \frac{1}{K}|Y|^{1.5}\mu_1$ with $K = (\frac{1}{K}|Y|^{1.5}\mu_1)[\R]$, $\mu_2$ is the law of $\exp(\Nc(-\frac12\sigma_1^2,\sigma_1^2))-1$ with $\sigma_1 = 0.1$, and $\nu_2$ is the law of $\exp(\Nc(-\frac12\sigma_2^2,\sigma_2^2))-1$ with $\sigma_2 = 0.2$. The scale indicates the mass in each point of the grid, the mass of the entropic approximation of the optimal coupling is the yellow zone. Notice that in all the cases the optimal coupling is supported on at most two maps. We saw this in all our experiment, we conjecture that for almost all $\mu,\nu$ this is the case.

Figure \ref{fig:left_curtain} shows well the left curtain coupling from \cite{beiglboeck2016problem} and \cite{henry2016explicit}. Figure \ref{fig:distance_cost} shows the optimal coupling for the distance cost. This coupling has been studied by Hobson \& Neuberger \cite{hobson2012robust}. They predict that this coupling is concentrated on two graphs. Finally, Figure \ref{fig:random_cost} shows how we may find solutions for any kind of cost function.

\subsubsection{Example in two dimensions}

\begin{figure}[h]
\centering
\begin{subfigure}{.24\textwidth}
  \centering
  \includegraphics[width=.99\linewidth]{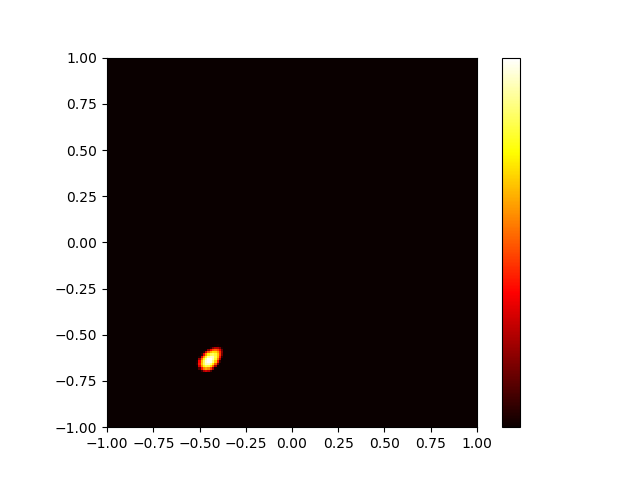}
  \caption{$X=-(0.45,0.65)$}
  \label{fig:sub1}
\end{subfigure}%
\begin{subfigure}{.24\textwidth}
  \centering
  \includegraphics[width=.99\linewidth]{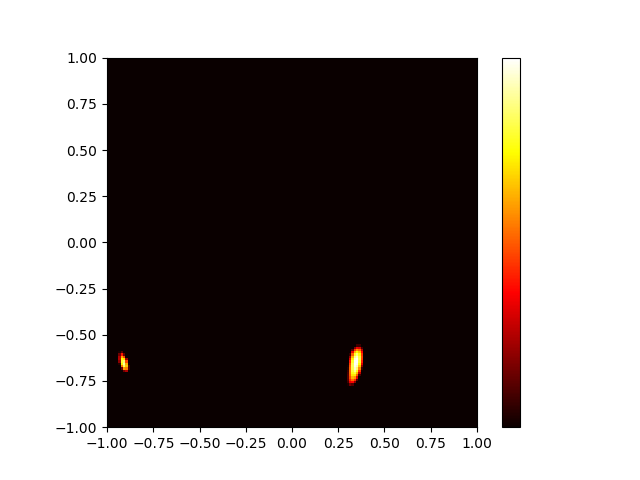}
  \caption{$X=(0.3,-0.66)$}
  \label{fig:sub2}
\end{subfigure}
\begin{subfigure}{.24\textwidth}
  \centering
  \includegraphics[width=.99\linewidth]{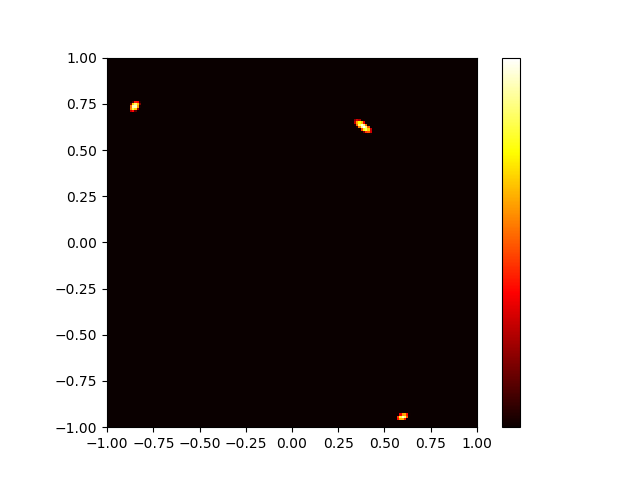}
  \caption{$X=(0.2,0.44)$}
  \label{fig:sub3}
\end{subfigure}
\begin{subfigure}{.24\textwidth}
  \centering
  \includegraphics[width=.99\linewidth]{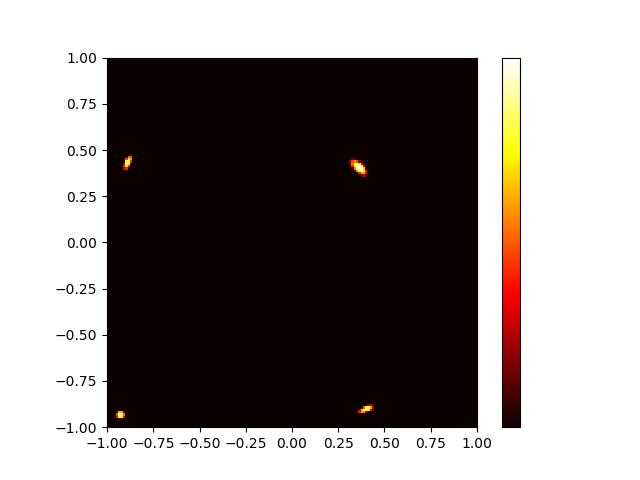}
  \caption{$X=(0.13,0.16)$}
  \label{fig:sub4}
\end{subfigure}
\caption{Optimal coupling conditioned to several values of $X$.}
\label{fig:test}
\end{figure}

In dimension $2$, it has been proved in \cite{de2018local} that for the cost function $c:(x,y)\in\R^2\x\R^2\longmapsto x_1(y_1^2+2y_2^2)+x_2(2y_1^2+y_2^2)$, the kernel of optimal probabilities are concentrated on the intersection of two ellipses with fixed characteristics, except for their position and their scale. Figure \ref{fig:test} is meant to test this theoretical result. We do an entropic approximation with a grid $160\x160$, and $\eps = 10^{-4}$. Then we selected $4$ points $x_1 := (-0.45,-0.65)$, $x_2 := (0.3,-0.66)$, $x_3 := (0.2,0.44)$, and $x_4 := (0.13,0.16)$ and draw the kernels of the approximated optimal transport $\P^*$ conditioned to $X = x_i$ for $i = 1,2,3,4$. We see on this figure that for $i = 1,2,3,4$, $\P^*(\cdot|X = x_i)$ is concentrated on exactly $i$ points, showing that all the numbers between $1$ and $4$ are reached. It seems that no trivial result can be proved on the number of maps that we may reach.

\bibliographystyle{plain}
\bibliography{bib/mabib}
\end{document}